\documentclass[twoside,bezier,reqno]{amsart}

\usepackage{epsfig}
\usepackage{amsmath,amsthm,amsfonts,amssymb,amscd,euscript}

\input{epsf}
\newcommand{\filebegin}{\begin{document}}
\newcommand{\fileend}{\end{document}}
\makeatother
\def\thefootnote{}
\newcommand{\lo}{\longrightarrow}
\newcommand{\NMM}{\hspace*{2mm}}

\renewcommand{\baselinestretch}{1.1}

\input{epsf}
\renewcommand{\baselinestretch}{1.1}
\def\n{\noindent}%

\numberwithin{equation}{section}
\def\mapdown#1{\Big\downarrow\rlap
{$\vcenter{\hbox{$\scriptstyle#1$}}$}}
\newtheorem{theorem}{Theorem}[section]
\newtheorem{lemma}[theorem]{Lemma}
\newtheorem{proposition}[theorem]{Proposition}
\newtheorem{corollary}[theorem]{Corollary}

\theoremstyle{definition}
\newtheorem{definition}[theorem]{Definition}
\newtheorem{example}[theorem]{\sc Example}
\newtheorem{xca}[theorem]{Exercise}
\theoremstyle{remark}
\newtheorem{remark}[theorem]{Remark}
\usepackage{tikz-cd}

\DeclareMathOperator{\sigmod}{\sigma-\left(R\text{-Mod}\right)}
\DeclareMathOperator{\sighh}{\left(\sigma\text{-}HH\right)}
\DeclareMathOperator{\sigbkn}{\sigma-\left(BKN\right)}
\DeclareMathOperator{\sigher}{\sigma\text{-}(R\text{-her})}
\DeclareMathOperator{\sigcoher}{\sigma\text{-}(R\text{-coher})}
\DeclareMathOperator{\signat}{\sigma\text{-}(R\text{-nat})}
\DeclareMathOperator{\sigconat}{\sigma\text{-}(R\text{-conat})}
\DeclareMathOperator{\sigTors}{\sigma\text{-}(R\text{-TORS})}
\DeclareMathOperator{\sigtors}{\sigma\text{-}(R\text{-tors})}

\usepackage[
breaklinks=true,
extension=pdf,
colorlinks=true,
linkcolor=black,
citecolor=red,
urlcolor=blue,
]{hyperref}

\begin{document}

\setcounter{page}{1} \noindent

\vspace*{2cm}
\begin{center}
{\bf\large On $\boldsymbol{\sigma}$-classes  of modules with applications}
 \\[0.5cm]
{Oscar Alberto Garrido-Jim\'enez$^a$, Hugo Alberto Rinc\'on-Mej\'ia$^{b*}$ \footnote{$^*$Corresponding Author} \\[2mm]
$^a$Facultad de Ciencias, Universidad Nacional Aut\'onoma de M\'exico Circuito Exterior, C.U., C.P. 04510, Ciudad de M\'exico, M\'exico\\
$^b$Facultad de Ciencias, Universidad Nacional Aut\'onoma de M\'exico Circuito Exterior, C.U., C.P. 04510, Ciudad de M\'exico, M\'exico\\[2mm]
{\tt E-mail: oagj@ciencias.unam.mx}\\
{\tt E-mail: hurincon@gmail.com}
} \\[2mm]
\end{center}%
\vspace*{0.5cm}
\begin{quotation}
\noindent
{\footnotesize
{\sc Abstract.}
In this paper we introduce some lattices of classes of left R-module relative to a preradical sigma. These lattices are generalizations of the lattices R-TORS, R-tors, R-nat, R-conat, of torsion theories, hereditary torsion theories, natural classes and conatural classes, respectively. We define the lattices $\sigTors$, $\sigtors$, $\signat$, $\sigconat$, which reduce to the lattices mentioned above, when one takes sigma as the  identity. We characterize the equality between these lattices by means of the $\sighh$ condition, which we introduce. We  also present some results about   $\sigma$-retractable rings,  $\sigma$-Max rings extending results about Mod-retractable rings and Max rings.  }
\end{quotation}

\n \textbf{2020 Mathematics subject classification: } 16D80, 
16S90, 
16S99, 
16W99. 

\markboth 
{Oscar Alberto Garrido-Jim\'enez, Hugo Alberto Rinc\'on-Mej\'ia}
 {On $\boldsymbol{\sigma}$-classes  of modules with applications}



\section{Introduction}

$R$ will denote an associative ring with 1 and $R$-Mod will denote the category of left $R$-modules and $R$-morphisms.
Several results of module theory are latticial in nature.  The behavior of the lattices associated with a ring determines several properties of the ring. The comparison between several lattices associated to a ring often provides interesting information. We will use this method with some lattices of module classes relative to a preradical $\sigma$.

Recall that a preradical in $R$-Mod is a functor $\sigma:R\text{-Mod}\to R\text{-Mod}$ such that for each $R$-module $M$ one has that $\sigma(M)\leq M$ and for each $R$-morphism $f:M\to N$ one has that $f\left(\sigma(M)\right)\leq \sigma(N). $ Note that a preradical $\sigma$ is a subfuntor of the identity functor in $R$-Mod. 
	Recall that preradicals n $R$-mod conforms the big complete lattice $R$-pr, where order is given  by $\sigma \leq \tau \Leftrightarrow \sigma(M)\subseteq \tau(M), \text{ for each }  _RM.$  A family of preradicals $\{r_i\}_I$     has a an infimum $ \wedge r_i$  which evaluated in a module $M$ gives $\bigcap_{i\in I} r_i (M)$, and   $\{r_i\}_I$  also has a supremun $ \vee r_i$  which evaluated in a module $M$  gives  $\sum_{i\in I} r_i (M)$. 
    The  identity in $R$-mod  functor is the greatest preradical and the preradical $\underline 0$  which assigns to each module its zero submodule is the least preradical.
    Recall that there are in $R$-pr two binary operations, composition, denoted by $\sigma\tau$  (defined by  $\sigma\tau(M)=\sigma\left(\tau(M)\right)$  and cocomposition, denoted by $\sigma:\tau$ (defined by  $(\sigma:\tau)(M)/\sigma(M)=\tau\left(M/\sigma(M)\right)$.

We say that a preradical $\sigma$ is  a left exact preradical if it   is left exact as a functor.  
$\sigma$ is idempotent if $\sigma \sigma =\sigma$.  $\sigma$ is stable (costable) if for each injective $R$-module $E$ (for each projective $R$-module $P$) one has that $E=\sigma(E)\oplus E'$  (one has that $P=\sigma(P)\oplus P'$) for some $E'\leq E$ (for some $P'\leq P$).
	
    There are two important classes of modules associated with a preradical $\sigma$, its torsion class $\mathbb{T }_\sigma=\{M \mid \sigma(M)=M\}$ and its torsion free class $\mathbb{F}_\sigma=\{M \mid \sigma(M)=0 \}.$
    $\mathbb{T }_\sigma$ is a class closed under the taking of coproducts and epimorphic images (i.e. it is a pretorsion class) and $\mathbb{F}_\sigma$  is a class closed under taking submodules and products (i.e is a pretorsion free class).
    
We say that an ordered pair $(\mathbb{T},\mathbb{F})$ of module classes is a torsion theory if $\mathbb{T}$ is closed under taking epimorphic images, coproducts and extensions, and $\mathbb{F}=\{N\mid Hom(M,N)=0,\ \forall M\in\mathbb{T}\}$ or equivalently, if $\mathbb{F}$ is closed under taking submodules, products, extensions and $\mathbb{T}=\{M\mid Hom(M,N)=0,\ \forall N\in\mathbb{F}\}$. 
A torsion theory $(\mathbb{T},\mathbb{F})$  having  $\mathbb{T}$   closed under taking submodules, or equivalently, having $\mathbb{F}$ closed under taking injective hulls is called a hereditary torsion theory.

The class of hereditary torsion theories is a complete lattice that has been extensively studied (see, for example \cite{golan} or \cite{stenstroem}).

For further information about lattices of classes of modules see \cite{ClassesofModules}.  

Comparisons between lattices of classes of modules may characterize types of rings. By example, in \cite{modretra}, the authors prove that  every torsion theory is hereditary, if and only if,  all modules are retractable. They called $R$-Mod-retractable ring to a ring with this property.
	A module $M$ is retractable when $Hom(M,N)\neq 0$ for each nonzero submodule $N$ of $M$.

In \cite{Associated to rings with respect to a preradical}, the lattice  $\sigTors$ of  $\sigma$-torsion classes, and the lattice  $\sigtors$ of  $\sigma$-hereditary torsion classes are defined for a preradical $\sigma$.
We refer to the reader to \cite[Chapter VI]{stenstroem} and \cite{BKN} for more information about  preradicals on $R$-Mod.

In \cite{inducidas por preradicales}, for an exact left and stable preradical $\sigma$,   the lattice $\signat$ of $\sigma$-natural classes is defined  and, for an exact and costable preradical $\sigma$,  the lattice $\sigconat$ of $\sigma$-conatural classes is defined.

We introduce, for a ring $R$ and a preradical $\sigma$, the condition $\sighh$.

We say that $R$ satisfies the condition $\sighh$ if for any pair of  $R$-modules $M, N$ with $\sigma(M)\neq 0$, and $\sigma(N)\neq 0$, one has that $Hom(\sigma(N),M)\neq 0$ if and only if $Hom(M,\sigma(N))\neq 0$. 
This condition is related to the parainjective $R$-modules  of $\sigma$-torsion and with the paraprojective $R$-modules  of $\sigma$-torsion. This condition is also important with respect to left $\sigma$-semiartinian rings, the left $\sigma$-max rings and the $\sigma$-retractable rings introduced in \cite{Associated to rings with respect to a preradical}.  
Finally, it turns out that the rings satisfying the $\sighh$-condition are the rings for which the lattices of $\sigma$-torsion classes, of  hereditary $\sigma$-torsion classes, of  $\sigma$-natural classes and of $\sigma$-conatural classes are all the same.

\section{Some kinds of preradicals}

Recall the following definitions.
\begin{definition}
A preradical  $\sigma$ is
\begin{enumerate}
\item  left exact if and only if $\sigma(A)=A\cap \sigma(B)$ whenever $A$ be a submodule of $B$ if and only if $\sigma$ is idempotent and its pretorsion class $\mathbb{T}_\sigma$ is hereditary.
\item  cohereditary if $\sigma$ preserves epimorphisms, this is equivalent to the pretorsion-free class of $\sigma$, $\mathbb{F}_\sigma$ being closed under quotients and  it is also equivalent to $\sigma(M)=\sigma(R)M$, for each $M.$
\item  splitting in $M$  if $\sigma(M)\leq_\oplus M$.
\item splitting if it splits in each module.
\item cosplitting if it is hereditary and cohereditary.
\item  stable if it splits in each injective module.
\item costable if it splits in each projective module.
\item centrally splitting if it has a complement $\sigma'$  in $R$-pr, or equivalently if there exists a central idempotent $e\in R$  such that $\sigma(M)=eM$, for each module $M$.
\end{enumerate}

\end{definition}
We include some basic properties of preradicals for the reader's convenience.

\begin{lemma} The following statements about a preradical $\sigma$ are equivalent. 
\begin{enumerate}
\item $\sigma$ is cohereditary.
\item $\sigma=\sigma (R)\cdot (-)$,  i.e. $\sigma(M)=\sigma (R)\cdot (M) $, for each module $M$.
\item $\sigma$ is a radical and it pretorsion free class $\mathbb{F_\sigma}$ is closed under quotients.
\end{enumerate}

\end{lemma}
\begin{proof}
$(1)\Rightarrow (2)$ Let $M$ be a module and let $f:R\rightarrow M$ be a morphism, $f=(-)\cdot x$ for some  $x\in M$. Then we have that the epimorphism $R\rightarrow Rx$ restricts to an epimorphism $\sigma(R)\rightarrow \sigma (Rx)$, thus $\sigma(R)x= \sigma(Rx)$. Now, take the epimorphism $R^{(M)}\twoheadrightarrow M$ induced by the family of morphisms $\{ (-)\cdot x : R\rightarrow M \}_{x\in M}$ then we have the epimorphism $\sigma(R^{(M)})\twoheadrightarrow \sigma(M)$. Thus $\sigma(M)=\sum_{x\in M} \sigma(R)x=\sigma(R)M$.

$(2)\Rightarrow (1)$ In general, if $I$ is an ideal of $R$ then $I\cdot(-)$ is a preradical preserving epimorphisms, as it is immediately verified. More over, $I\left(M/IM\right)=0$ shows that $I\cdot(-)$ is a radical.

$(2) \Rightarrow (3)$ Assume that $f:M \twoheadrightarrow N$ is an epimorphism with \(M\in \mathbb{F}_\sigma\) then we have an epimorphism $0=\sigma(M)\twoheadrightarrow \sigma(N)$. $\sigma$ is a radical as it was noted in the above argument.

$(3)\Rightarrow(1)$ Let $f:M\twoheadrightarrow N$ be an epimorphism.  $f$ induces an epimorphism $\overline{f}:M/\sigma(M)\twoheadrightarrow N/f\left(\sigma(M)\right).$ Thus we have the  situation of the following diagram.

\[
\begin{tikzcd}
0 \arrow[r] & f\left(\sigma(M)\right)\arrow[r, tail] & N\rar[two heads] & N/f\left(\sigma(M)\right)\arrow[r]& 0 \\
0\arrow[r]& \sigma(M) \rar[u, two heads]{f[} \rar[tail]{} & M\rar[u, two heads]{f} \arrow[r, two heads]& M/\sigma(M)\rar[u, two heads]{\overline{f}} \arrow[r]& 0.
\end{tikzcd}
\]

As $\sigma$ is a radical, it follows that  $\sigma\left(M/\sigma(M)\right)=0$. Now, as $\mathbb{F}_\sigma$ is a class closed under taking quotients, it follows that $\sigma\left(N/f\left(\sigma(M)\right)\right)=0$. Then, as $f\left(\sigma(M)\right)\leq N$and as  $\sigma$ is a radical, it follows that  $\sigma\left(N/f\left(\sigma(M)\right)\right)=\sigma(N)/f\left(\sigma(M)\right)$ (see,  \cite[Lemma 1.1, Chap. VI]{stenstroem}). So, $\sigma(N)/f\left(\sigma(M)\right)=0$ and hence $\sigma(M)=f\left(\sigma(M)\right).$ Therefore  $\sigma$ is cohereditary.
 \end{proof}
 
\begin{lemma} \label{her and coher sigma R}
A preradical $\sigma$ is hereditary and cohereditary if and only if $\sigma = \sigma(R)\cdot(-)$  and the ideal $ \sigma(R)$  satisfies that $x \in \sigma(R)x $, for each $x\in\sigma(R)$.
\end{lemma}
\begin{proof}
$\Rightarrow]$ We have seen that $\sigma$ cohereditary implies  $\sigma = \sigma(R)\cdot(-)$.  On the other hand, as $\sigma$ is hereditary then it is idempotent thus $\sigma(R)$ is idempotent and $\mathbb{T}_\sigma$ is hereditary class.  So if $x\in\sigma(R)$ then $Rx\leq \sigma(R)\in\mathbb{T}_\sigma$, applying $\sigma$,  $x\in  Rx=\sigma(Rx)=\sigma(R)Rx=\sigma(R)x.$

$\Leftarrow]$ Conversely, if $\sigma$ is a cohereditary radical  and   $x \in \sigma(R)x $, for each $x\in\sigma(R)$, then $\sigma(R)\leq\sigma(R)\sigma(R)\leq \sigma(R).$ Thus $\sigma$  is an idempotent radical. We only have to show that $\mathbb{T}_\sigma$ is hereditary class, which is equivalent to $\mathbb{F}_\sigma$ being closed under injective hulls. Suppose that \(M\in \mathbb{F}_\sigma\) but \(E(M)\notin \mathbb{F}_\sigma \) . Then \(\sigma(E(M))\neq 0\) and then  \(M\cap\sigma(R)(E(M)) \neq 0\).  Then there exists \(0\neq m\in \sigma(R)(E(M)) \), then \(m=a_1x_1+\cdots+a_nx_n\), with \(a_i\in \sigma(R)\)  and \(x_i \in E(M)\). Let us take one such \(m\) where \(n\) be  least possible.  Notice that \(n\neq1\), because if  \(0\neq m=a_1x_1\) then \(a_1=b_1a_1\)  for some \(b_1 \in \sigma (R)\), by hypothesis. So, \(0\neq m =b_1a_1x_1=b_1m=0\),  a contradiction. Now,  in \(m=a_1x_1+\cdots+a_nx_n\), we can write \(a_i=b_ia_i\)   for some \(b_i\in\sigma(R)\). Note that \(m=b_1a_1x_1+\cdots+b_na_nx_n\) and that \(0=b_1m=b_1a_1x_1+\cdots+b_1a_nx_n\) thus making a subtraction we have that \(m=(b_1-b_2)a_2x_2+\cdots+(b_1-n_n)a_nx_n \)  with less than \(n\) summands, contradicting the choice of \(m\) and \(n\). This proves that \(\mathbb{F}_\sigma\) is closed under taking injective hulls, and this completes the proof.
\end{proof}

\begin{lemma} A preradical $\sigma$ is exact, i.e. $\sigma$ preserves exact sequences of modules if and only if it is hereditary and cohereditary.
\end{lemma}
\begin{proof}
If $\sigma$ is exact, clearly is left exact and right exact, thus it  is hereditary and cohereditary.
If $\sigma$ is cohereditary and hereditary then it preserves epimorphisms and it is left exact. Thus  \(\sigma\) preserves short exact sequences and it is well known that this implies that $\sigma$ is an exact functor (See \cite[Chap. I, \textsection 5, Lemma 5.1]{stenstroem}).
\end{proof}

\begin{lemma} \label{esrsdcs}  For a preradical $\sigma$ the following conditions are equivalent:
\begin{enumerate}
\item  $\sigma$ is exact and and $\sigma(R)$  is a direct summand of $R$.
\item $\sigma$ centrally splits.
\end{enumerate}
\end{lemma}

\begin{proof}
1) $\Rightarrow 2)$ As \(\sigma\)  is exact then \(\sigma= \sigma(R)\cdot (-)\) and  $\sigma (R)$ is an idempotent ideal having the property that \(x\in \sigma (R) x, \forall x \in \sigma (R)\). Let us call \(\rho\) the preradical which assigns to a module \(M\) the submodule \(\{x\in M \mid \sigma (R)x =0\}\).  \(\rho\)  is a left exact radical whose torsion class \(\mathbb{T}_\rho\) is closed under taking products.  \(\rho(R)\) is a two sided ideal.  
As \(R=\sigma(R)\oplus J\) for some left ideal \(J\), by hypothesis, and as \(\sigma(R)J \subseteq \sigma(R)\cap J =0\)  then, \(J \subseteq \rho (R)\). Thus \(R =\sigma(R)\oplus \rho (R) \). Then $ \sigma(R)=Re $  with \(e\)  a central idempotent. Thus \(\sigma= \sigma(R)\cdot (-)=e\cdot (-)\). Hence \(\sigma\) centrally splits.

2)$ \Rightarrow$ 1)  This follows from Lemma \ref{her and coher sigma R} 
\end{proof}

 \begin{lemma}
 The following statements are equivalent

\begin{enumerate}
\item  $\sigma$ is exact and stable.
\item $\sigma$ centrally splits.
\item $\sigma$ is exact and costable.
\end{enumerate}
\end{lemma}

\begin{proof}
Is clear that if  \(\sigma\)  centrally splits, then  it is  exact. Now, as \(\sigma(M)\)  is a direct summand of \(M\), for each module \(M\), then \(\sigma\) is stable and costable.

Conversely, if \(\sigma\) is exact and stable, then \(\sigma(E(R))\hookrightarrow  E(R)\)  splits, thus  \(\sigma(R)\hookrightarrow  R\) also splits. And if \(\sigma\) is exact and costable, then \(\sigma(R)\hookrightarrow  R\)  splits, because \(R\) is a projective module. We conclude using Lemma \ref{esrsdcs}. 
\end{proof}

\bigskip

It is a known fact that the class of left exact preradicals, which we will denote $R$-lep (by left exact preradicals) is in bijective correspondence with the set of linear filters for the ring and it is also in bijective correspondence with the class of Wisbauer classes.
 A Wisbauer class is a class of modules closed under taking submodules, quotients and direct sums.  
Thus $R$-lep is a complete lattice.
We denote by $R$-pr the big lattice of preradicals in $R$-mod.

\section{$\sigma$-classes of modules}

In \cite{inducidas por preradicales}, the authors introduced  lattices of module classes induced by a preradical 
$\sigma$. Namely they introduced $\sigma$-hereditary classes, $\sigma$-cohereditary classes, $\sigma$-natural classes and $\sigma$-conatural classes and showed that they conform (big) lattices, which were denoted  $\sigma\text{-}(R\text{-her})$, $\sigma\text{-}(R\text{-coher})$, $\sigma\text{-}(R\text{-nat})$ and $\sigma\text{-}(R\text{-conat})$, respectively. We recall the definitions and some relevant results  for the sequel.

\begin{definition}
Let  $\sigma\in R$-pr. A class $\mathcal{C}$   of $R$-modules is a \textbf{ $\boldsymbol{\sigma}$-hereditary class} if it satisfies the following conditions:
\begin{enumerate}
    \item $\mathbb{F}_\sigma\subseteq\mathcal{C}$.
    \item For each $M\in\mathcal{C}$ and each $N\leq M$ it follows that $\sigma(N)\in\mathcal{C}.$    
\end{enumerate}
We denote by  \textbf{$\boldsymbol{\sigher}$} to the  collection of all   $\sigma$-hereditary classes.
\end{definition}

Following Golan, we call skeleton of a lattice $\mathcal{L}$ with a smallest element, to the set of pseudocomplements of $\mathcal{L}$ and denote it Skel($\mathcal{L}$).

\begin{lemma}{\cite[Lemma 1]{inducidas por preradicales}} Let $\sigma\in R$-pr. The collection  $\sigher$ is a  pseudocomplemented big lattice. Moreover, if $\sigma$ is idempotent, then $\sigher$ is strongly pseudocomplemented and for each  $\mathcal{C}\in \sigher$ its pseudocomplement, denoted by $\mathcal{C}^{\perp_{\leq_\sigma}},$ is given by 
\begin{align*}
    \mathcal{C}^{\perp_{\leq_\sigma}}& =\left\{M\in R\text{-Mod}\ |\ \forall\ N\leq M,\ \sigma(N)\in\mathcal{C}\Rightarrow N\in\mathbb{F}_\sigma\right\}.
\end{align*}
In addition, $Skel\left(\sigher\right)$ is a boolean lattice.
\end{lemma}
\bigskip

Given the above lemma, we can make the following definition.

\begin{definition}
Let $R$ be a ring and let  $\sigma\in R$-pr be left exact and stable. We define the lattice of $\sigma$-natural classes,  denoted by \textbf{$\boldsymbol{\signat}$},  as
\[\sigma\text{-}\left(R\text{-nat}\right)=Skel(\sigher).\]
\end{definition}

\begin{definition}
Let $\sigma\in R$-pr. A class $\mathcal{C}$ of $R$-modules  is a \textbf{ $\boldsymbol{\sigma}$-cohereditary class} if it satisfies the following conditions:
\begin{enumerate}
    \item $\mathbb{F}_\sigma\subseteq\mathcal{C}$.
    \item For each $M\in\mathcal{C}$ and each epimorphism $M\twoheadrightarrow L$ we have that $\sigma(L)\in\mathcal{C}.$    
\end{enumerate}
We denote the collection of all   $\sigma$-cohereditary classes by \textbf{$\boldsymbol{\sigcoher}$}.
\end{definition}

The following proposition is a generalization of \cite[Proposition 3.6]{seudocomplements}.

\begin{proposition}{\cite[Proposition 11]{inducidas por preradicales}} Let $\sigma\in R$-pr. If $\sigma$ is idempotent and cohereditary, then $\sigcoher$ is a strongly pseudocomplemented big lattice  and for each  $\mathcal{C}\in \sigcoher$ the pseudocomplement, of $\mathcal{C}$, denoted by $\mathcal{C}^{\perp_{/_\sigma}},$ is given by 
\begin{align*}
    \mathcal{C}^{\perp_{/_\sigma}}=\left\{M\in R\text{-Mod}\ |\ \forall\ M\twoheadrightarrow L,\ \sigma(L)\neq 0\Rightarrow\sigma(L)\notin\mathcal{A}\right\}\cup\mathbb{F}_\sigma.
\end{align*}
In addition, $Skel\left(\sigcoher\right)$ is a boolean lattice. 
\end{proposition}

\begin{definition}
Let $R$ be a ring and let  $\sigma\in R$-pr be exact and costable. We define the lattice of $\sigma$-conatural classes, denoted by  \textbf{$\boldsymbol{\sigconat}$} as
\[\sigconat=Skel(\sigcoher).\]
\end{definition}

Later, in \cite{Associated to rings with respect to a preradical}, the authors introduced the lattices $\sigTors$ and $\sigtors$ as follows:

\begin{definition}
Let $\sigma\in R$-pr. A class $\mathcal{C}$ of $R$-modules    \textbf{$\boldsymbol{\sigma}$-torsion  class} if $\mathcal{C}$ is a $\sigma$-cohereditary class, closed under taking arbitrary direct sums and under taking extensions,  if it is also  a  $\sigma$-hereditary class, then we will say that  $\mathcal{C}$ is an  \textbf{hereditary $\sigma$-torsion class}.
\end{definition}

\begin{definition}
Let $\sigma\in R$-pr. We denote the collection of all  $\sigma$-torsion classes by \textbf{$\boldsymbol{\sigTors}$} and by \textbf{$\boldsymbol{\sigtors}$} the collection of all  hereditary $\sigma$-torsion classes. 
\end{definition}

In \cite{Associated to rings with respect to a preradical}, the authors define the following assignations relative to a preradical $\sigma$:
\begin{enumerate}
    \item $\sigma^*:\mathcal{P} (R\mbox{-Mod})\to\mathcal{P} (R\mbox{-Mod})$ given by $\sigma^*(\mathcal{C})=\{\sigma(M)\ |\ M\in\mathcal{C}\}$ and
    \item $\overleftarrow{\sigma}:\mathcal{P} (R\mbox{-Mod})\to\mathcal{P} (R\mbox{-Mod})$ given $\overleftarrow{\sigma}(\mathcal{C})=\{M\ |\ \sigma(M)\in\mathcal{C}\}.$
\end{enumerate}
In addition, it holds that $\sigma^*\left(\overleftarrow{\sigma}\left(\sigma^*\left(\mathcal{C}\right)\right)\right)=\sigma^*(\mathcal{C})$ and  $\overleftarrow{\sigma}\left(\sigma^*\left(\overleftarrow{\sigma}\left(\mathcal{C}\right)\right)\right)=\overleftarrow{\sigma}(\mathcal{C}),$ for each $\mathcal{C}\in \mathcal{P} (R\mbox{-Mod}).$

The authors showed the following results.
\begin{proposition}\cite[Corollary 3.10]{Associated to rings with respect to a preradical} If  $\sigma$ is a radical, then \[\sigtors=\{\mathcal{C}\in R\text{-tors }|\ \mathbb{F}_\sigma\subseteq\mathcal{C}\}.\]
\end{proposition}

\begin{theorem}\cite[Theorem 3.13]{Associated to rings with respect to a preradical} If $\sigma$ is an exact preradical, then \[\sigTors=\{\overleftarrow{\sigma}(\mathcal{C})\ |\ \mathcal{C}\in R\text{-TORS }\}.\]
\end{theorem}

$\sigma$-retractable modules and $\sigmod$-retractable rings also introduced in \cite{Associated to rings with respect to a preradical}.

\begin{definition}

Let $R$  be a ring and $\sigma \in  R\text{-}pr$. We say that an $R$-module $M$   is
\textbf{$\boldsymbol{\sigma}$-retractable} if for every submodule $N \leq  M$ with $\sigma(N) \neq 0$ we have that $Hom (M, \sigma (N)) \neq 0$ 
and we will say that $R$ is a \textbf{$\boldsymbol{\sigma\text{-}(R\text{-}Mod)}$-retractable ring} if every $R$-module is $\sigma$-retractable.
\end{definition}

In \cite{On retractability and its}, the  authors introduced the condition $(HH)$ and characterized  the rings satisfying it, via the  retractability and properties of lattices of module classes, namely, the lattice of natural classes, the lattice of conatural classes, and the lattice of hereditary torsion theories. In the following, we extend some of these concepts and some of these results.

\begin{definition}
Let $R$ be a ring and $\sigma  \in  R\text{-}pr$. We will say that $R$ satisfies the \textbf{condition $\boldsymbol{\sighh}$}
if for any modules $M, N$ with $\sigma(M)\neq 0,\  \sigma(N) \neq 0,$ it holds that
$Hom(\sigma (N), M) \neq  0 \Leftrightarrow Hom(M, \sigma(N)) \neq  0.$ 
\end{definition}

\begin{remark}
Note that if  $\sigma=1_{R\text{-}Mod}$, then we have condition $(HH)$ (see \cite{On retractability and its}),  i.e. for any  $R$-modules $M$ and $N$ we have that   
\[Hom(N,M)\neq 0\Leftrightarrow Hom(M,N)\neq 0.\]
\end{remark}

Recall that a module $M$ is a  \textbf{parainjective module}   if for each module   $N$, whenever there is some monomorphism  $f\in \mathrm{Hom}(M,N)$ there exists an epimorphism  $g\in \mathrm{Hom}(N,M)$. A module $M$ is a \textbf{paraprojective module}  if for each module $N$, whenever there is some   epimorphism $f\in \mathrm{Hom}(N,M)$ there exists a monomorphism  $g\in \mathrm{Hom}(M,N)$.

 $R\text{-}Simp$ will denote a family of  representatives of isomorphism classes of  left simple  modules.

\begin{remark} \label{consecuen de SHH}Let $R$ be a ring and $\sigma \in R$-pr. If $R$ satisfies the $\sighh,$ condition, then:

\begin{enumerate}
\item\label{SHH impli Sretra}  $R$ is a $\sigmod$-retractable ring.
\item\label{SHH impli parainyec}  For any $S\in R\text{-}Simp\cap\mathbb{T}_\sigma$, it follows that $S$ is parainjective.
\item\label{SHH impli parapro}  For any $S\in R\text{-}Simp\cap\mathbb{T}_\sigma$, it follows that \(S\) is paraprojective. \end{enumerate}
\end{remark}
\begin{proof}
\ \begin{enumerate}
\item  It is clear.
\item Let $S\in R\text{-}Simp \cap\mathbb{T}_\sigma$ and let  \(M\) be an \( R\)-module such that \(S\) embeds in \(M\).  Since \( S \in  \mathbb{T}_{\sigma}, \) it follows that \(\sigma(S) = S\). Since $S \in \mathbb{T}_\sigma$,  then 
\(\sigma(S) = S\), so $(Hom(\sigma(S), M) \neq  0$. Then, $Hom(M, \sigma(S)) \neq 0,$ since \(R\) satisfies the $\sigma\text{-}HH$ condition. Thus, $Hom(M, S) \neq 0,$ i.e., \(S\) is a quotient of \(M\), whence S is parainjective.
\item Let $S\in R\text{-}Simp \cap\mathbb{T}_\sigma$ and let  \(M\) be an \(R\)-module such that \(S\) is a quotient of \(M\). It follows that
\(Hom(M, \sigma(S)) \neq 0\), since \(S = \sigma(S)\) because \( S \in \mathbb{T}_{\sigma}\). Then, from the $\sighh$ condition, it follows that $Hom(\sigma(S),M)\neq 0$ , i.e.  $Hom(S, M) \neq  0,$ whence \(S\) embeds  in \(M\). Therefore, \( S\) is paraprojective.
 \end{enumerate}\ \ \ \ \ \ \ \ \ \ \ \ \ \ \ \ \ \ \ \ \ \ \ \ \ \ \ \ \ \ \ \ \ \ \ \ \ \ \ \ \ \ \ \ \ \ \ \ \ \ \ \ \ \ \ \ \ \ \ \ \ \ \ \ \ \ \ \ \ \ \ \ \ \ \ \ \ \ \ \ \ \ \ \ \ \ \ \ \ \ \ \ \ \ \ \ \ \ \ \ \ \ \ \ \ \ \ \ \ \ \ \ \ \ \ \ \ \ \ \ \ \ \ \ \ \ \ \ \ \ \ \ \end{proof}

\begin{definition}

  Let $R$ be a ring and $\sigma \in  R\text{-}lep.$ We will say that an \(R\)-module $M \neq  0$  is $\sigma$-atomic if any nonzero submodule $N$ of $M$ with $\sigma(N) \neq 0$ contains a submodule
$S\in R\text{-}Simp \cap\mathbb{T}_\sigma,$ and we will say that $R$ is a left $\sigma$-semiartinian ring if any \(R\)-module is $\sigma$-atomic.\end{definition}
Recall the following definition which were introduced in \cite{Associated to rings with respect to a preradical}.

\begin{definition}
Let $R$ be a ring and take an idempotent preradical  $\sigma$.  We say that a non zero $R$-module $M $ is $\sigma$-coatomic if for each  non zero quotient $L$ of $M$ with $\sigma(L)\neq 0$, $L$ has a quotient $S$,  $S\in  R\text{-}Simp\cap\mathbb{T}_\sigma$. We will say that  $R$ is a left $\sigma$-max ring, if each  $R$-module is  $\sigma$-coatomic.
\end{definition}


\begin{theorem}\label{Sparain Smax}
Let $R$ be a ring and let $\sigma$ be a left exact preradical.  If each $S$, such that $S\in  R\text{-}Simp\cap\mathbb{T}_\sigma$ is parainjective, then $R$ is a left $\sigma$-max  ring.
\end{theorem}
\begin{proof}
Let $M\neq 0$ be an $R$-module and  $L$ be a non zero quotient of $M,$ with $\sigma(L)\neq 0.$ Let us consider  a cyclic  submodule $Rx$ of  $\sigma(L)$. As  $\sigma$ is left exact, then  $\sigma$ is idempotent and $\mathbb{T}_\sigma$ is a hereditary class. So  $\sigma(L)\in\mathbb{T}_\sigma,$ and thus  $Rx\in\mathbb{T}_\sigma.$  Now, if we  consider a simple quotient of  $Rx,$  $S$, say, we have that  $S\in\mathbb{T}_\sigma,$ thus $S\in  R\text{-}Simp\cap\mathbb{T}_\sigma$. 
Taking the  inclusion of $S$ in its injective hull  $E(S),$ we obtain a non zero  morphism from $Rx$ to $E(S).$ This morphism has a non zero  extension  $f:L\to E(S).$  Note that  $S\leq  f(L)$  because $S$ is essential in $E(S).$ As  $S$ is parainjective, then it is a quotient of $f(L)$ and thus $S$  is a quotient of  $L.$ See the below diagram.
\[
\begin{tikzcd}
M \arrow[r, twoheadrightarrow] & L \ar[rrdd, dashed, "f"] & &  \\
& \sigma(L) \arrow[u, tail] & & \\
& Rx \arrow[u, tail] \arrow[r, two heads] & S \arrow[r, tail] & E(S).
\end{tikzcd}
\]

\end{proof}

\begin{theorem}\label{Sparapro Ssemi}
Let  $R$ be a ring and let  $\sigma$ be a left exact preradical. If each  simple module  $S$  in  $\mathbb{T}_\sigma$  is paraprojective, then  $R$ is a left $\sigma$-semiartinian ring.
\end{theorem}
\begin{proof}
Let $M$ be an  $R$-module and let $N$ be a non zero submodule of $M$ with  $\sigma(N)\neq 0$. Let us take a non zero cyclic submodule  $Rx$ of $\sigma(N).$ As  $\sigma$ is left exact, then  $\sigma(N)\in\mathbb{T}_\sigma$, thus  $Rx\in\mathbb{T}_\sigma.$ 

\[
\begin{tikzcd}
\sigma(N) \arrow[r, tail] & N \arrow[r, tail] & M \\
Rx \arrow[u, tail] \arrow[r, two heads] & S &
\end{tikzcd}
\]

Now, if $S$ is a simple quotient of  $Rx,$ then we have that  $S\in\mathbb{T}_\sigma$.  By hypothesis, $S$ is paraprojective, thus  $S$ embeds in $Rx$ and from this,  $S$ embeds in $M.$
\end{proof}

\begin{corollary}\label{SHH Smax y Ssemi}
Let  $R$ be a ring and let  $\sigma\in R$ be a left exact preradical. If $R$ is a ring satisfying the $\sighh$-condition,
then $R$ is a left  $\sigma$-max and a left  $\sigma$-semiartinian ring.
\end{corollary}
\begin{proof}
It follows from  \eqref{SHH impli parainyec} and \eqref{SHH impli parapro} of Remark \ref{consecuen de SHH} and from  Theorems \ref{Sparain Smax} and  \ref{Sparapro Ssemi}.
\end{proof}

\begin{theorem}\label{sigHH sii pariny y parapro}
Let  $R$ be a ring and let $\sigma$ be a left exact preradical. The following statements are equivalent:
\begin{enumerate}
    \item $R$ has the $\sighh$-condition.
    \item $R$ is a  $\sigmod$-retractable ring and each $S\in  R\text{-}Simp\cap\mathbb{T}_\sigma$ is paraprojective.
    \item Any $S\in  R\text{-}Simp\cap\mathbb{T}_\sigma$ is paraprojective and parainjective.
\end{enumerate}
\end{theorem}
\begin{proof}
$(1)\Rightarrow (2)$ It follows  from \eqref{SHH impli Sretra} and \eqref{SHH impli parapro}  of   Remark \ref{consecuen de SHH}.\\
\noindent $(2)\Rightarrow (3)$ This is clear.\\
\noindent $(3)\Rightarrow (1)$ Let  $M$ and $N$  be two  $R$-modules such that  $\sigma(M),\sigma(N)\neq 0$ and let us suppose that $ f\in Hom\left( \sigma(N),M \right)\neq 0$. We have the following diagram

\[
\begin{tikzcd}
\sigma(N) \rar[two heads]{f \downharpoonright} & f\left(\sigma(N)\right) \arrow[r, tail] & M \arrow[dl, dashed, two heads] \\
f^{-1}(S) \arrow[u, tail] \rar[ two heads]{f\downharpoonright} & S \arrow[u, tail] &
\end{tikzcd}
\]
As  $\sigma$ is left exact, we have that  $\sigma(N)\in\mathbb{T}_\sigma$, and thus  $0\neq f\left(\sigma(N)\right)\in\mathbb{T}_\sigma$. Now, from Theorem \ref{Sparapro Ssemi}, there exists  $S\in  R\text{-}Simp\cap\mathbb{T}_\sigma,$ a  submodule of $f\left(\sigma(N)\right)$, which is  parainjective. Thus $S$ is a quotient of  $M$. 
Besides, $S$ is a quotient of  $f^{-1}(S)$, hence   $S$ embeds in $f^{-1}(S)$, because  $S$ is paraprojective. Thus we have that  $Hom\left(M,\sigma(N)\right)\neq 0.$

Let us now see that if $Hom\left(M,\sigma(N)\right)\neq 0,$ then  $Hom\left(\sigma(N),M\right)\neq 0.$ To do so, let us consider  $0\neq f\in Hom\left(M,\sigma(N)\right).$ Given that  $\sigma\in R$-lep and $0\neq f(M)\leq\sigma(N),$ then  $f(M)\in\mathbb{T}_\sigma,$ i.e., $\sigma\left(f(M)\right)=f(M)\neq 0.$ Thus, by Theorem \ref{Sparapro Ssemi}, there exists $S\in R\text{-Simp}\cap\mathbb{T}_\sigma$ such that  $S$ embeds in $\sigma\left(f(M)\right).$

\[
\begin{tikzcd}
M \rar[two heads]{f \downharpoonright} & f\left(M\right) \arrow[r, tail] & \sigma(N) \arrow[dl, dashed, two heads] \\
f^{-1}(S) \arrow[u, tail] \rar[ two heads]{f\downharpoonright} & S \arrow[u, tail] &
\end{tikzcd}
\]
Then, taking into account that every  $\sigma$-torsion simple module is parainjective, we have that is a quotient of $\sigma(N).$ Finally, using the fact that every $\sigma$-torsion simple module  is paraprojective, we have that  $S$ embeds in $f^{-1}(S)$,Thus we obtain a nonzero morphism from $\sigma(N)$ to $M.$ 
\end{proof}
\bigskip

A ring  $R$ is called a  \textbf{$\boldsymbol{BKN}$-ring} if between any two nonzero modules there exists a nonzero morphism. Such rings are introduced and characterized  in \cite{BKN}. In the following, we generalize these rings.

\begin{definition}
Let $R$  be a ring and let  $\sigma\in R$-pr.  We say that  $R$ is a \textbf{$\boldsymbol{\sigbkn}$-ring} if for any  $R$-modules  $M$ and $N$ with   $\sigma(M),\sigma(N)\neq 0,$ it holds that  
\[Hom\left(\sigma(N),\sigma(M)\right)\neq 0.\]
\end{definition}

Recall that a ring $R$ is called \textbf{left  local ring} if there is only one isomorphism type of simple $R$-modules. In the following definition we generalize this concept.

\begin{definition}
Let $R$ be a ring and  $\sigma\in R$-pr. We say that a ring  $R$ is  \textbf{left $\boldsymbol{\sigma}$-local ring}  if\[\left| R\text{-Simp}\cap\mathbb{T}_\sigma\right|=1.\]
\end{definition}

\begin{theorem}
Let $R$ be a ring and let  $\sigma$ be an exact and costable preradical. The following statements are equivalent:
\begin{enumerate}
    \item $R$ is a $\sigbkn$-ring.
    \item $R$ satisfies the $\sighh$-condition and it is a left  $\sigma$-local ring.
\end{enumerate}
\end{theorem}
\begin{proof}
$(1)\Rightarrow(2)$ Let  $M$ and $N$ be two  $R$-modules such that  $\sigma(M),\sigma(N)\neq 0$ and with  $Hom\left(\sigma(N),M\right)\neq 0.$ As $R$ is a $\sigbkn$-ring, we have that $Hom(\sigma(M),\sigma(N))\neq 0$. As  $\sigma$ is exact and costable, then  $\sigma(M)$ is a quotient of $M.$ Hence  $Hom(M,\sigma(N))\neq 0.$
\bigskip

Now, if  $S_1,S_2\in R\text{-Simp}\cap\mathbb{T}_\sigma,$ we have that   $Hom\left(\sigma(S_1),\sigma(S_2)\right)\neq 0,$ thus  $Hom(S_1,S_2)\neq 0$. Then   $S_1 \cong S_2.$   Hence  $R$    is a left $\sigma$-local ring .
\bigskip

\noindent $(2)\Rightarrow(1)$ Let  $M$ and $N$ be two  $R$-modules such that  $\sigma(M),\sigma(N)\neq 0.$ From Corollary   \ref{SHH Smax y Ssemi}, we have that  $R$ is $\sigma$-semiartinian, so there exist  $S_1,S_2\in R\text{-}Simp\cap\mathbb{T}_\sigma$ such that  $S_1$ embeds in $\sigma(M)$ and $S_2$ embeds in $\sigma(N).$ As  $R$ satisfies  the   $\sighh$ condition, it follows that  $Hom\left(\sigma(M),S_1\right)\neq 0.$ Finally, as  $R$  is a left  $\sigma$-local ring,  we see that $S_1\cong S_2$ and from this, we get  $Hom\left(\sigma(M),\sigma(N)\right)\neq 0.$ 
\end{proof}
\bigskip

The following definition appears in \cite[Definition 22]{on natural}.
\begin{definition}
A module class  $\mathcal{C}$ satisfies the condition \emph{($CN$)} if whenever 
    each  nonzero quotient of an arbitrary module $M$ shares a nonzero quotient with some element of $\mathcal{C}$, it happens that $M$  belongs to $\mathcal{C}$.
\end{definition}

In the following definition we  generalize the $(CN)$ condition. 

\begin{definition}\label{condición sigma-CN}
Assume  $\mathcal{C}\subseteq R$-Mod.  $\mathcal{C}$ satisfies the \textbf{$\boldsymbol{(\sigma\text{-}CN)}$-condition} if
\[\left(\begin{array}{c}
\text{for each epimorphism }\ M\twoheadrightarrow L\ \text{with}\ \sigma(L)\neq 0,\\ \text{there exist}\ C\in\mathcal{C},
N\in R\text{-Mod},\ \text{with}\ \sigma(N)\neq 0,\\ \text{and epimorphims }\ L\twoheadrightarrow N\twoheadleftarrow C
\end{array}\right)\Longrightarrow M\in\mathcal{C}.\]
\end{definition}

\begin{lemma}\label{cocientes de elementos de C están en C}
Let $\mathcal{C}$ be the class of $R$-modules  satisfying the  $(\sigma\text{-}CN)$-condition and  let  $C\in\mathcal{C}.$ If  $L$ is a quotient of  $C$ with  $\sigma(L)\neq 0$, then  $L\in\mathcal{C}.$
\end{lemma}
\begin{proof}
Let $H$ be a quotient of $L$ with  $\sigma(H)\neq 0.$ Note that  $H$ is  also a quotient of $C$. As  $\mathcal{C}$ satisfies  $(\sigma\text{-}CN)$, it follows that   $L\in\mathcal{C}.$
\end{proof}
\bigskip

The following theorem is a generalization of  \cite[Theorem 22]{on natural}.
\begin{theorem}
Let  $\sigma$-pr be an exact and  costable preradical and  $\mathcal{C}\subseteq R$-Mod. The following statements are equivalent for $\mathcal{C}$:
\begin{enumerate}
    \item $\mathcal{C}\in\sigconat$.
    \item $\mathcal{C}$ satisfies  $(\sigma\text{-}CN).$
    \item $\mathcal{C}\in\mathcal{L}_{/_\sigma}$ and  $\mathcal{C}=\left(\mathcal{C}^{\perp_{/_\sigma}}\right)^{\perp_{/_\sigma}}.$
\end{enumerate}
\end{theorem}
\begin{proof}
$(1)\Rightarrow (2)$ Let $\mathcal{C}\in\sigma\text{-}\left(R\text{-conat}\right)$ and  $M$ be an $R$-module. As $\mathcal{C}\in\sigma\text{-}\left(R\text{-conat}\right)$, then $\mathcal{C}=\mathcal{A}^{\perp_{/_\sigma}},$ for some $\mathcal{A}\in \mathcal{L}_{/_\sigma},$ thus by \cite[Proposition 11]{inducidas por preradicales}, we have
\[\mathcal{A}^{\perp_{/_\sigma}}=\left\{M\in R\text{-Mod}\ |\ \forall\ M\twoheadrightarrow L,\ \sigma(L)\neq 0\Rightarrow\sigma(L)\notin\mathcal{A}\right\}\cup\mathbb{F}_\sigma.\]
Note that if for each non zero quotient  $L$ of $M$, with  $\sigma(L)\neq 0$ it holds that $\sigma(L)\notin\mathcal{A},$  then $M\in\mathcal{A}^{\perp_{/_\sigma}},$  $M\in\mathcal{C}.$  Let us assume that  $L$ is a non zero quotient of $M$, with $\sigma(L)\neq 0,$ and such that  $\sigma(L)\in\mathcal{A}$ and that there exist  $C\in\mathcal{C}$ and  $N\in R$-Mod, with $\sigma(N)\neq 0,$ with epimorphisms  $L\twoheadrightarrow N\twoheadleftarrow C.$ As  $C\in\mathcal{C}$ and  $C\twoheadrightarrow N$, is an epimorphism with $\sigma(N)\neq 0,$ we get that  $\sigma(N)\notin\mathcal{A}.$ 

Besides for the epimorphism $f:L\twoheadrightarrow N$, as  $\sigma$ is exact, then $\sigma$ is idempotent and cohereditary, thus  $\sigma(N)=\sigma\left(\sigma(N)\right)$ and    $f_{\downharpoonright} :\sigma(L)\twoheadrightarrow\sigma(N)$ is also an epimorphism. As $\mathcal{A}\in\mathcal{L}_{/_\sigma}$ it follows that  $\sigma(N)\in\mathcal{A}$, a contradiction. Thus  $\sigma(L)\notin\mathcal{A}$ and hence $M\in\mathcal{C}.$

$(2)\Rightarrow(3)$ Let $\mathcal{C}$  be a  class of $R$-modules satisfying  $(\sigma\text{-}CN).$ Firstly, we show that $\mathcal{C}\in\mathcal{L}_{/_\sigma}.$ For this, let us consider  $C\in\mathcal{C}$ and let  $L\in R$-Mod be a quotient of  $M$ with  $\sigma(L)\neq 0.$ As  $\sigma$ is  an exact costable preradical, we have that  
\begin{align}\label{sigma de C es un cociente de C}
    C=\sigma(R)C\oplus C'=\sigma(C)\oplus C',
\end{align}
where $C'=\{c\in C\ |\ \sigma(R)c=0\}$. It follows that  $\sigma(C)$ is a quotient of $C$. Then, by Lemma \ref{cocientes de elementos de C están en C}, $\sigma(L)\in\mathcal{C}.$

Now, to show that  $\mathcal{C}=\left(\mathcal{C}^{\perp_{/_\sigma}}\right)^{\perp_{/_\sigma}},$ let us recall the following  description of the double pseudocomplement of $\mathcal{C},$
\begin{equation}\label{descripción del doble seudocomplemento}
\begin{split}
    \left(\mathcal{C}^{\perp_{/_\sigma}}\right)^{\perp_{/_\sigma}}=\{M\in R\text{-Mod}\ |\ &\forall M\twoheadrightarrow K,\ \sigma(K)\neq 0  \\& \Rightarrow\exists \sigma(K)\twoheadrightarrow L,\ \sigma(L)\neq 0\ \text{and}\ \sigma(L)\in\mathcal{C}\}\cup\mathbb{F}_\sigma,
\end{split}
\end{equation}
see \cite[Definition 4]{inducidas por preradicales}.
So, if  $C\in\mathcal{C}$ and $K$ is a  quotient of $C$ with $\sigma(K)\neq 0,$ then as  $\sigma$  is idempotent and cohereditary, we get an epimorphism 
\begin{align*}
    \sigma(C)\twoheadrightarrow\sigma(K).
\end{align*}
In the same way that in  \eqref{sigma de C es un cociente de C}, $\sigma(C)$ is a quotient of $C,$ thus from Lemma \ref{cocientes de elementos de C están en C} $\sigma(C)\in\mathcal{C}.$ Applying  Lemma \ref{cocientes de elementos de C están en C} one more time, we obtain that $\sigma(K)\in\mathcal{C}.$
From  \eqref{descripción del doble seudocomplemento} we get that  $C\in\left(\mathcal{C}^{\perp_{/_\sigma}}\right)^{\perp_{/_\sigma}}.$ This shows that $\mathcal{C}\subseteq\left(\mathcal{C}^{\perp_{/_\sigma}}\right)^{\perp_{/_\sigma}}.$

For the other inclusion, let us  consider $M\in\left(\mathcal{C}^{\perp_{/_\sigma}}\right)^{\perp_{/_\sigma}}$ and a quotient $L$ of   $M$ with  $\sigma(L)\neq 0.$ As  $\sigma(L)$ is a  quotient of $L$, then $\sigma(L)$ is a quotient of $M.$ Also $\sigma\left(\sigma(L)\right)=\sigma(L)\neq 0$. So as $M\in\left(\mathcal{C}^{\perp_{/_\sigma}}\right)^{\perp_{/_\sigma}},$ there exists a quotient $T$ of  $\sigma(L)$ with $\sigma(T)\neq 0$ and $\sigma(T)\in\mathcal{C}.$ Finally,  as  $\mathcal{C}$ satisfies $(\sigma\text{-}CN),$ we get that  $M\in\mathcal{C}.$ So, $\left(\mathcal{C}^{\perp_{/_\sigma}}\right)^{\perp_{/_\sigma}}\subseteq\mathcal{C},$ hence  $\mathcal{C}=\left(\mathcal{C}^{\perp_{/_\sigma}}\right)^{\perp_{/_\sigma}}.$

$(3)\Rightarrow(1)$ This is  clear.
\end{proof}

\begin{remark}\label{descripcion de la sigmaconatural generada}
Let  $\sigma$ be an  exact and  costable preradical and let  $\mathcal{C}\subseteq R$-Mod. If  $\xi_{\sigma \text{-conat}}\left(\mathcal{C}\right)$ denotes the least $\sigma$-conatural  class containing  $\mathcal{C}$, then
\begin{equation*}
    \begin{split}
        \xi_{\sigma\text{-conat}}\left(\mathcal{C}\right)=& \{M\in R\text{-Mod}\ |\ \forall M\twoheadrightarrow L,\ \sigma(L)\neq 0  \ \text{there exist }\ C\in\mathcal{C}, \\ & N\in R\text{-Mod}, \text{with}\ \sigma(N)\neq 0,\ \text{such that }\ L\twoheadrightarrow N\twoheadleftarrow C \}\cup \mathbb{F}_\sigma
    \end{split}
\end{equation*}
\end{remark}

In \cite[Theorem 42]{on natural} left max rings are characterized via conatural classes.  The following theorem generalizes this result.

\begin{proposition}\label{en sigmamax conat esta determinada por simples}
Let  $\sigma$ be an exact and  costable preradical. The following statements are equivalent:
\begin{enumerate}
    \item $R$ is a left  $\sigma$-max ring.
    \item Each  $\sigma$-conatural class is generated by a family of simple  $\sigma$-torsion modules.
\end{enumerate}  
\end{proposition}
\begin{proof}
$(1)\Rightarrow(2)$ Let  $\mathcal{C}$ be a non trivial  $\sigma$-conatural class and  let $\mathcal{S}\subseteq\mathcal{C}$ be the class of all simple   $\sigma$-torsion modules  in $\mathcal{C}.$ Notice that  $\mathcal{S}$ is non empty because  $R$ is a $\sigma$-max ring. Now, if $0\neq M\in\mathcal{C}$ and  $L$ is a non zero quotient of $M$  with $\sigma(L)\neq 0,$ as  $R$ is  $\sigma$-max, then   $L$ has a simple quotient $S\in R\text{-Simp}\cap \mathbb{T}_\sigma.$ By Lemma \ref{cocientes de elementos de C están en C}, $S\in\mathcal{C}$, by  Remark \ref{descripcion de la sigmaconatural generada}, we have that  $M\in\xi_{\sigma\text{-conat}}(\mathcal{S}).$  It follows that $\mathcal{C}=\xi_{\sigma\text{-conat}}(\mathcal{S}).$

$(2)\Rightarrow(1)$ Let  $M\in R$-Mod and  $L$  be a non zero quotient of $M$ with $\sigma(L)\neq 0.$ By hypothesis, there exists $\mathcal{S}\subseteq R\text{-Simp}\cap \mathbb{T}_\sigma$ such that  $\xi_{\sigma\text{-conat}}(M)=\xi_{\sigma\text{-conat}}(\mathcal{S}).$ As $L\in \xi_{\sigma\text{-conat}}(M),$ $L$ has a quotient  $S\in\mathcal{S}$, Thus  $R$ is a left  $\sigma$-max ring.
\end{proof}
\bigskip

The following theorem generalizes \cite[Proposition 2.11]{on boolean}.

\begin{theorem}\label{sigmaconat contenida en sigmators}
Let  $\sigma$ be an exact and  costable preradical. The following statements are equivalent:
\begin{enumerate}
    \item Each $S\in R\text{-Simp}\cap \mathbb{T}_\sigma$ is parainjective.
    \item $\sigconat\subseteq \sigtors$.
    \item $\sigconat\subseteq \mathcal{L}_{\leq_\sigma}.$
\end{enumerate}
\end{theorem}
\begin{proof}
$(1)\Rightarrow(2)$ In \cite[Proposition 6.5]{Associated to rings with respect to a preradical} it is shown that for an exact and costable  preradical  $\sigma,$ and for a left  $\sigma$-max ring, each   $\sigma$-conatural class is closed under taking direct sums. Now, by hypothesis,  each  $S\in R\text{-Simp}\cap \mathbb{T}_\sigma$ is parainjective, from  Theorem \ref{Sparain Smax}, we obtain that  $R$ is a left  $\sigma$-max ring. Thus it suffices to show that any  $\sigma$-conatural class is a $\sigma$-hereditary class.

For this, let us consider  $\mathcal{C}\in \sigconat$, $M\in\mathcal{C}$ and $N\leq M$ with $\sigma(N)\neq 0.$ Let  $S\in\text{Simp}\cap\mathbb{T}_\sigma$ be a quotient of $\sigma(N),$ which exists because  $R$ is a left $\sigma$-max ring. If  $T$  be the push-out of  $\sigma(N)\twoheadrightarrow S$ and  $\sigma(N)\hookrightarrow M,$ as in the following diagram
\[\begin{tikzcd}
\sigma(N) \arrow[d, twoheadrightarrow] \arrow[r,rightarrowtail] & N \arrow[r,rightarrowtail] & M \arrow[d,dashed, twoheadrightarrow]\\
S \arrow[rr, tail, dashed] & & T,
\end{tikzcd}\]
then $S$ embeds in  $T.$ Thus, as  $S$ is parainjective by hypothesis, then  $S$ is a simple  $\sigma$-torsion quotient of  $T.$ Then $S$ is also a simple  $\sigma$-torsion quotient of $M.$ Then any simple $\sigma$-torsion quotient of  $\sigma(N)$ belongs to $\mathcal{C}$.  Thus, by  Proposition \ref{en sigmamax conat esta determinada por simples}, $\sigma(N)\in\mathcal{C}.$ Hence $\mathcal{C}$ is a  $\sigma$-hereditary class.

$(2)\Rightarrow(3)$ This is clear.

$(3)\Rightarrow(1)$ Let  $S\in\text{Simp}\cap\mathbb{T}_\sigma$ and let  $M$  be an  $R$-module such that  $S$ embeds in $M.$ By hypothesis, we have that  $S\in\xi_{\sigma\text{-conat}}(M).$ From Remark \ref{descripcion de la sigmaconatural generada} it follows that  $S$ is a quotient  of $M,$ then  $S$ is parainjective.
\end{proof}

\begin{remark}
If  $\sigma$ is an exact and costable  preradical,   then for each  $R$-module $M$, we have a decomposition $M=\sigma(M)\oplus M',$ where  $M'=\{m\in M\ |\ \sigma(R)m=0\}.$ Notice than $\sigma$  has a complement $\sigma'$, in $R$-pr satisfying  
\[\sigma'(M)=M',\]
for  each $M\in R$-Mod.
\end{remark}

\begin{theorem}
Let  $\sigma$ be an exact and  costable preradical. If  $R$ satisfies $\sighh$ and  $\left(\sigma'\text{-}HH\right),$ then $R$ satisfies  $(HH).$ 
\end{theorem}
\begin{proof}
Let $M$ and $N$ be two  $R$-modules and suppose that  $Hom(N,M)\neq 0.$  Let us  take a non zero $f\in Hom(N,M)$. If  $\nu:M=\sigma(M)\oplus\sigma'(M)\to\sigma(M)$ denotes the natural projection on  $\sigma(M)$ and  $\nu f\neq 0,$ then $Hom(N,\sigma(M))\neq 0.$ So, by  $\sighh, $there exists a non zero morphism  $ g:\sigma(M)\to N,$ thus  $g\nu:M\to N$ is a non zero morphism, $Hom(M,N)\neq 0.$

Now, if $\nu f=0,$ then taking the co-restriction of $f$ to its image, we have that  $ f\upharpoonright :N\to\sigma'(M)$ is a non zero morphism. From condition  $\left(\sigma'\text{-}HH\right),$ there exists a non zero morphism  $h:\sigma'(M)\to N.$ Finally, $h\nu' :M\to N$ is  non zero morphism, where  $\nu':M\to\sigma'(M)$ is the natural  projection over  $\sigma'(M).$ Thus, $Hom(M,N)\neq 0.$\end{proof}
\bigskip

In \cite[Theorem 1.13]{On retractability and its} it is shown that if  $R$ satisfies condition $(HH),$ then each torsion class is a conatural class. The following corollary is a consequence.
\begin{corollary}\label{sigma y sigma' TORS-CONAT}
Let  $\sigma\in R$-pr be exact and costable. If  $R$ satisfies condition  $\sighh$ and  condition $\left(\sigma'\text{-}HH\right),$ then each torsion class is a conatural class.
\end{corollary}

\begin{theorem}
\label{sigmaTORS contenida en sigmaconat con sigma y sigma'}
Let $\sigma\in R$-pr be exact and  costable. If  $R$ satisfies  condition $\sighh$ and condition $\left(\sigma'\text{-}HH\right)$, then \[\sigTors\subseteq \sigconat.\]
\end{theorem}
\begin{proof}
Let   $\mathcal{C}\in \sigTors.$ As $\sigma$ is  exact, there exists $\mathcal{D}\in R$-TORS such that  $\mathcal{C}=\overleftarrow{\sigma}(\mathcal{D})$, see \cite[Theorem 3.13]{Associated to rings with respect to a preradical}.
Now, from  Corollary  \ref{sigma y sigma' TORS-CONAT}, we have that $\mathcal{D}\in R$-conat. Then, as  $\sigma$ is an exact  and costable preradical, we have that  $\overleftarrow{\sigma}(\mathcal{D})\in R\text{-}(\sigma\text{-conat}),$ see \cite[Proposition 15]{inducidas por preradicales}, i.e., $\mathcal{C}\in \sigconat$. 
\end{proof}

\begin{theorem}
Let  $\sigma\in R$-pr be exact and costable. If  $R$ satisfies condition $\sighh$ and $\left(\sigma'\text{-}HH\right)$, then  \[\sigTors= \sigconat=\sigtors.\]
\end{theorem}
\begin{proof}
It follows from  Theorem \ref{sigmaTORS contenida en sigmaconat con sigma y sigma'}, from  Remark \ref{consecuen de SHH}, from  Theorem \ref{sigmaconat contenida en sigmators} and from the fact  that $\sigtors\subseteq \sigTors.$
\end{proof}

\begin{lemma}\label{Submodulo de sig-torsion en C}
Let  $\sigma$ be an  exact and  costable preradical, let $R$ be a ring with  $(\sigma\text{-}HH)$, $\mathcal{C}\in \sigTors$ and  $M\in R$-Mod. If  $M$ has a quotient  $L$ with  $L\in \mathcal{C}$ and  $\sigma(L)\neq 0,$ then $M$  has a  non zero submodule  $N,$ such that  $N\in\mathcal{C}\cap\mathbb{T}_\sigma.$ Moreover, there exist a largest submodule  $K$  of $M$ with the property that  $K\in\mathcal{C}\cap\mathbb{T}_\sigma.$
\end{lemma}
\begin{proof}
As $\sigma$  is an exact and costable preradical, then $\sigma(L)$ is a quotient of $L.$ Now, as $\mathcal{C}$ is a  $\sigma$-torsion class and $\sigma(L)\neq 0,$  it follows that $\sigma(L)\in\mathcal{C}.$ 

As $\sigma(L)$ is a nonzero quotient of $M,$ then $Hom(M,\sigma(L))\neq 0.$ Thus, from the $(\sigma\text{-}HH)$-condition, there exists an morphism $0\neq h:\sigma(L)\to M.$ Let $N=h\left(\sigma(L)\right),$ and note that $N$ is a nonzero submodule of $M$. Moreover, as $N$ is a quotient of $\sigma(L),$ then $\sigma(N)$ is a quotient of $\sigma\left(\sigma(L)\right),$ because $\sigma$ is cohereditary, as $\sigma\left(\sigma(L)\right)=\sigma(L)$, then $\sigma(N)=N$ with $\sigma(N)\in\mathcal{C}.$ Hence, $N\in\mathcal{C}\cap\mathbb{T}_\sigma.$

Now, as $\mathcal{C}$ is a class closed under taking direct sums, then $\mathcal{C}_M=\{ N\leq M\ |\ 0\neq N\in \mathcal{C}\cap\mathbb{T}_\sigma\}\neq\emptyset$ and as $\sigma$ is left exact, we get that 
\[\bigoplus_{C\in\mathcal{C}_M}C\in\mathcal{C}\cap\mathbb{T}_\sigma.\]
As $\sum_{C\in\mathcal{C}_M}C$ is a quotient of $\bigoplus_{C\in\mathcal{C}_M}C$ and as $\sigma$ is a cohereditary preradical, it follows that  $\sum_{C\in\mathcal{C}_M}C\in\mathcal{C}\cap\mathbb{T}_\sigma.$ Hence, $K=\sum_{C\in\mathcal{C}_M}C$ is the largest $\sigma$-torsion submodule of $M$ belonging to $\mathcal{C}.$ 
\end{proof}

\begin{lemma}\label{Tors ent sig* sigconat}
Let $\sigma$ be an  exact and  costable preradical, let  $R$  be a ring  with the $\sighh$-condition and let  $\mathcal{C}\subseteq R$-Mod. If $\mathcal{C}\in \sigTors$, then $\sigma^*(\mathcal{C})\in R\text{-conat}$.
\end{lemma}
\begin{proof}
Let us show that $\sigma^*(C)$ has $(CN).$ For this, let us take an $R$-module $M$ such that for any non zero quotient $N$ there exist $C\in\mathcal{C}$ and $L\in R$-Mod as in the following diagram: 
\[
\begin{tikzcd}
M\arrow[r, two heads] & N \arrow[r, two heads] & L & \sigma(C)\arrow[l, two heads] 
\end{tikzcd}
\]
and we will see that $M\in\sigma^*(\mathcal{C}).$

In this situation, $\sigma(C)\in\mathcal{C}$, because $\sigma(C)$ is a quotient of $C$, as $\sigma$ is exact and costable, and $\mathcal{C}$ is a cohereditary class. Now, note that, $L=\sigma(L)$ because $\sigma$ is cohereditary and $\sigma(C)=\sigma\left(\sigma(C)\right),$ thus $L\in\mathcal{C}.$ So, by Lemma \ref{Submodulo de sig-torsion en C}, we consider the largest submodule, $K$ of $M$ such that $K\in\mathcal{C}\cap\mathbb{T}_\sigma.$ 
Note that $K\neq 0 $  and that if $K=M,$ then we are done, because $M=\sigma(K)$ with $K\in\mathcal{C},$ that is, $M\in\sigma^*(\mathcal{C}).$
Assume to the contrary that $K\lneq M,$  then $M/K\neq_{R}0.$  By hypothesis,  $M/K$ also has a non zero submodule  $U\in\mathcal{C}\cap\mathbb{T}_{\sigma},$  as in the following commutative diagram

\[
\begin{tikzcd}
0 \arrow[r] & K\arrow[r, tail] & M\rar[two heads]{\pi} & M/K\arrow[r]& 0 \\
0\arrow[r]& K \uar[equal] \rar[tail]{} & \pi^{-1}\left(U\right)\arrow[u] \arrow[r, two heads]& U\arrow[u] \arrow[r]& 0.
\end{tikzcd}
\]

But as $\mathcal{C}\cap\mathbb{T}_{\sigma}$  is a class closed under extensions, then $\pi^{-1}\left(U\right)$ is a submodule of $M$  in $\mathcal{C}\cap\mathbb{T}_{\sigma},$  properly containing $K$, a contradiction. Then $K=M.$
\end{proof}
\bigskip

Let us recall that for any class $\mathcal{D}\subseteq R$-Mod, we have that $\overleftarrow{\sigma}(\mathcal{D})=\overleftarrow{\sigma}\sigma^*\overleftarrow{\sigma}(\mathcal{D})$; and that if $\sigma$ is an exact preradical, then a class $\mathcal{C}$ is $\sigma$-torsion if and only if $\mathcal{C}=\overleftarrow{\sigma}(\mathcal{D})$, for a torsion class $\mathcal{D}$. Besides, if $\sigma$ is an exact and costable preradical and $\mathcal{D}$ is a conatural class, then $\overleftarrow{\sigma}(\mathcal{D})$ is a $\sigma$-conatural class, see \cite{inducidas por preradicales} and \cite{Associated to rings with respect to a preradical}.


\begin{theorem}
\label{sigmaTORS contenida en sigmaconat}
Let $\sigma $ be exact and costable. If $R$ satisfies the $\sighh$ condition, then \[\sigTors\subseteq \sigconat.\]
\end{theorem}
\begin{proof}
Let $\mathcal{C}\in \sigTors.$ There exists a torsion class $\mathcal{D}$ such that $\mathcal{C}=\overleftarrow{\sigma}(\mathcal{D}).$ Thus, we have,

\[\mathcal{C}=\overleftarrow{\sigma}(\mathcal{D})=\overleftarrow{\sigma}\sigma^*\overleftarrow{\sigma}(\mathcal{D}).\]
It suffices to show that $\sigma^*\overleftarrow{\sigma}(\mathcal{D})$ is a conatural class. To see this, by Lemma \ref{Tors ent sig* sigconat} it suffices to show that $\overleftarrow{\sigma}(\mathcal{D})$ is a $\sigma$-torsion class. But this holds  by hypothesis.
\end{proof}

\begin{theorem}\label{sigHH ent sigTors=sigtors=sigconat}
Let $\sigma$ be an exact and costable preradical. If $R$ satisfies the $\sighh$ condition, then \[\sigTors= \sigconat=\sigtors.\]
\end{theorem}
\begin{proof}
It follows from Theorem \ref{sigmaconat contenida en sigmators}, from Theorem \ref{sigmaTORS contenida en sigmaconat} and the fact that $\sigtors\subseteq \sigTors.$
\end{proof}
\bigskip

The  following definition is dual to  Definition \ref{condición sigma-CN}.

\begin{definition}
Let  $\mathcal{C}\subseteq R$-Mod. We will say that  $\mathcal{C}$ satisfies  \textbf{condition $(\boldsymbol{\sigma\text{-}N)}$} if 
\[\left(\begin{array}{c}
\text{for each monomorphism}\ L\rightarrowtail M\ \text{with}\ \sigma(L)\neq 0,\\ \text{there exist}\ C\in\mathcal{C},
N\in R\text{-Mod},\ \text{with}\ \sigma(N)\neq 0\\ \text{and monomorphisms}\ L\leftarrowtail N\rightarrowtail C
\end{array}\right)\Longrightarrow M\in\mathcal{C}.\]
\end{definition}

\begin{lemma}\label{submodulos de elementos de C están en C}
Let $\mathcal{C}$ be a class of  $R$-modules satisfying the $(\sigma\text{-}N)$ condition and let  $C\in\mathcal{C}.$ If $L$ embeds in  $C$, with  $\sigma(L)\neq 0$, then $L\in\mathcal{C}.$
\end{lemma}
\begin{proof}
 Let $K$ be an  $R$-module embedded in  $L,$ with $\sigma(K)\neq 0.$ Then  $K$ embeds in  $C$ and as  $C\in\mathcal{C},$ we conclude that  $K\in\mathcal{C}.$
\end{proof} 

\begin{remark}
In \cite[Corollary 2]{inducidas por preradicales} it is shown that  
\[\sigma\text{-}\left(R\text{-nat}\right)=\mathcal{L}_{\{\leq_\sigma,\oplus,\sigma(E( )),ext\}}.\]
\end{remark}

\begin{theorem}
Let $\mathcal{C}\subseteq R$-Mod and let  $\sigma\in R$-pr be a left exact and stable preradical.  Then  $\mathcal{C}\in\sigma\text{-}\left(R\text{-nat}\right)$  if and only if  $\mathcal{C}$ satisfies  $(\sigma\text{-}N).$
\end{theorem}

\begin{proof}
Suppose that  $\mathcal{C}\subseteq R$-Mod is a $\sigma$-natural class. By definition of  $\sigma\text{-}\left(R\text{-}nat\right)$, one has that  $\mathcal{C}=\mathcal{A}^{\perp_{\leq_\sigma}},$ for some class  $\mathcal{A}\in\mathcal{L}_{\leq_\sigma}.$   Then, by \cite[Lemma 1]{inducidas por preradicales}, we have that
\[\mathcal{A}^{\perp^{/_\sigma}}=\left\{M\in R\text{-Mod}\ |\ \forall\ L\rightarrowtail M,\ \sigma(L)\in\mathcal{A}\Rightarrow L\in\mathbb{F}_\sigma\right\}.\]
Let $M$ be an $R$-module. Note that, if for any $L\rightarrowtail M$ it holds that  $\sigma(L)\in\mathcal{A}\Rightarrow L\in\mathbb{F}_\sigma,$ then $M\in\mathcal{A}^{\perp_{\leq_\sigma}}=\mathcal{C}.$ Suppose then that there exists $L\rightarrowtail M$ with  $\sigma(L)\in\mathcal{A},$ but that  $L\notin\mathbb{F}_\sigma.$ Additionally, suppose that there exists  $N\in R$-Mod, with  $\sigma(N)\neq 0$  and  $C\in\mathcal{C}$ such that  $C\leftarrowtail N\rightarrowtail L.$ 

Note that  $\sigma(N)\in\mathcal{C},$ since  $\mathcal{C}\in\mathcal{L}_{\leq_\sigma}$ and $C\in\mathcal{C}.$ It follows that  $\sigma(N)\notin \mathcal{A},$ since otherwise $N\in\mathbb{F}_\sigma.$

On the other hand, since  $\sigma$  is left exact, and hence it is idempotent, it follows that. $\sigma(N)=\sigma\left(\sigma(N)\right)\rightarrowtail\sigma\left(\sigma(L)\right)=\sigma(L).$ Finally, since $\mathcal{A}\in\mathcal{L}_{\leq_\sigma}$ and $\sigma(L)\in\mathcal{A}$  it follows that  $\sigma(N)\in\mathcal{A},$ which is a contradiction. Thus, there does not exist $L\rightarrowtail M$ with  $\sigma(L)\in\mathcal{A}$ and such that  $L\notin\mathbb{F}_\sigma.$ Therefore $M\in\mathcal{C}.$
\bigskip

Suppose now that $\mathcal{C}\subseteq R$-Mod is a class satisfying condition $(\sigma\text{-}N)$ and
let us see that, indeed,  $\mathcal{C}$ is a $\sigma$-natural class.

Let us begin by showing that $\mathcal{C}\in\mathcal{L}_{\leq_\sigma}.$  Let     $M\in\mathcal{C}$   and   $N\leq M.$  Note that if  $\sigma(N)=0,$ then  $N\in\mathcal{C},$ since $\mathbb{F}_\sigma\subseteq\mathcal{C}.$  Suppose then that  $\sigma(N)\neq 0$  and let  $L\rightarrowtail \sigma(N),$ with  $\sigma(L)\neq 0.$  Since  $L\rightarrowtail M$, there exist $H\in R$-Mod, with $\sigma(H)\neq 0,$ and  $C\in\mathcal{C}$ such that  $L\leftarrowtail H\rightarrowtail C.$

Now let us show that  $\mathcal{C}\in\mathcal{L}_{\oplus}.$ Let  $\{M_i\}_{i\in I}\subseteq\mathcal{C}$ and $L\rightarrowtail \oplus_{i\in I}M_i,$ with  $\sigma(L)\neq 0.$ Since  $\sigma$ is left exact, one has that $\sigma(L)\rightarrowtail\sigma\left(\oplus_{i\in I}M_i\right)=\oplus_{i\in I}\sigma(M_i).$ Now, by the projection argument, we have that there exist  $l\in \sigma(L)$ and $0\neq m_i\in \sigma(M_i)$, for some  $i\in I,$ such that  $Rl\cong Rm_i.$ Note also that $\sigma$ is idempotent and  $\mathbb{T}_\sigma$ is a hereditary class, since  $\sigma$ is left exact. It follows that  $\sigma(Rl)=Rl\neq 0,$ since $Rl\cong Rm_i\leq \sigma(M_i)$ and  $\sigma(M_i)\in\mathbb{T}_{\sigma}.$ Now, since  $M_i\in\mathcal{C}$ and  $Rm_i\rightarrowtail M_i,$ with  $\sigma(Rm_i)\neq 0,$ then there exist $C\in\mathcal{C}$  and  $N\in R$-Mod, with  $\sigma(N)\neq 0,$ such that  

\[
\begin{tikzcd}
L & Rl\cong Rm_i \arrow[l, tail] & N\arrow[l, tail]\arrow[r, tail]& C. 
\end{tikzcd}
\]
That is, $\oplus_{i\in I}M_i\in\mathcal{C}.$

Now, let us see that $\mathcal{C} \in \mathcal{L}_{\sigma (E())}.$  Let $M\in\mathcal{C}$  and  $L\rightarrowtail\sigma\left(E(M)\right).$ Since  $\sigma$ is left exact and stable, one has that  $\sigma\left(E(M)\right)=E\left(\sigma(M)\right).$  Now, given that  $\sigma(M)\leq_e E\left(\sigma(M)\right)$ and  $\sigma(L)\neq 0,$ it follows that  $\sigma(L)\cap\sigma(M)\neq 0.$  But, note that  $\sigma\left(\sigma(L)\cap\sigma(M)\right)=\sigma\left(E(M)\right)\cap\sigma(L)\cap\sigma(M)=\sigma(L)\cap\sigma(M)\neq 0,$ since  $\sigma$  is left exact. Therefore, exist $C\in\mathcal{C}$ and  $N\in R$-Mod,  with  $\sigma(N)\neq 0,$ such that 
\[
\begin{tikzcd}
L & \sigma(L)\cap\sigma(M) \arrow[l, tail] & N\arrow[l, tail]\arrow[r, tail]& C. 
\end{tikzcd}
\]
Whereupon, $\sigma\left(E(M)\right)\in\mathcal{C}.$

Finally, let us see that $\mathcal{C}\in\mathcal{L}_{ext}.$ Let 
\[
\begin{tikzcd}
0 \arrow[r] & M'\rar[tail]{f} & M\rar[two heads]{g} & M''\arrow[r]& 0
\end{tikzcd}
\]
be an exact sequence, with  $M',M''\in\mathcal{C}$, and $L\rightarrowtail M$ with  $\sigma(L)\neq 0.$ Without loss of generality, suppose that  $f$ is the  inclusion and that  $L\leq M.$ Note that, if  $\sigma(L)\cap M'\neq 0,$ then  $\sigma\left(\sigma(L)\cap M'\right)=\sigma(L)\cap\sigma(L)\cap M'=\sigma(L)\cap M'\neq 0,$ since  $\sigma$ is left exact. In this case, we have that there exists $C\in\mathcal{C}$ and $N\in R$-Mod, such that  $L\leftarrowtail\sigma(L)\cap M'\leftarrowtail N\rightarrowtail C,$ since  $M'\in\mathcal{C}.$

Now, if $\sigma(L)\cap M'=0,$ then $g\rfloor:\sigma(L)\rightarrowtail M''$ is a monomorphism, whence , there exist $C\in\mathcal{C}$ and $N\in R$-Mod, such that $L\leftarrowtail\sigma(L)\cong g\left(\sigma(L)\right)\leftarrowtail N\rightarrowtail C,$ since $M''\in\mathcal{C}.$
\end{proof}

\begin{remark}\label{descripcion de la sigmanatural generada}
Let  $\sigma\in R$-pr be left exact and stable, and  $\mathcal{C}\subseteq R$-Mod. If  $\xi_{\sigma\text{-nat}}\left(\mathcal{C}\right)$ denotes   the smallest  class $\sigma$-natural containing  $\mathcal{C}$, then
\begin{equation*}
    \begin{split}
        \xi_{\sigma\text{-nat}}\left(\mathcal{C}\right)=\{M\in R\text{-Mod}\ |\ &\forall L\rightarrowtail M,\ \sigma(L)\neq 0  \ \text{there exist}\ C\in\mathcal{C}, N\in R\text{-Mod}, \\ & \text{con}\ \sigma(N)\neq 0,\ \text{such that}\ L\leftarrowtail N\rightarrowtail C \}\cup \mathbb{F}_\sigma
    \end{split}
\end{equation*}
\end{remark}

\begin{proposition}\label{en sigmasemiar nat esta determinada por simples}
Let  $\sigma\in R$-pr be left exact and stable. The following statements are equivalent:
\begin{enumerate}
    \item $R$ is a left $\sigma$-semiartinian ring.
    \item Each   $\sigma$-natural class is generated by  a  family $\sigma$-torsion  simple modules.
\end{enumerate}
\end{proposition}
\begin{proof}
$(1)\Rightarrow(2)$ Let  $\mathcal{C}$ be a non trivial  $\sigma$-natural class and  $\mathcal{S}\subseteq\mathcal{C}$ be the class of all simple modules of $\sigma$-torsion of  $\mathcal{C}.$ The  class $\mathcal{S}$ is  non empty since  $R$ is a left  $\sigma$-semiartinian ring. Consider $ 0\neq M\in\mathcal{C} $, and   $_R L\neq 0$ embedding in $M$, with $\sigma(L)\neq 0.$ Since  $R$  is a left  $\sigma$-semiartinian ring, there exists $S\in R\text{-Simp}\cap \mathbb{T}_\sigma$ a submodule of  $L.$ By Lemma \ref{submodulos de elementos de C están en C}, the submodules of modules in  $\mathcal{C}$  belong to   $\mathcal{C}$,  so $S\in\mathcal{C}.$ Then, by  Remark \ref{descripcion de la sigmanatural generada}, one has that $M\in\xi_{\sigma\text{-nat}}(\mathcal{S}).$  It follows that $\mathcal{C}=\xi_{\sigma\text{-nat}}(\mathcal{S}).$

$(2)\Rightarrow(1)$ Let $M\in R$-Mod and $L$ be a nonzero submodule of  $M$ with $\sigma(L)\neq 0.$ By  hypothesis one has that  there exists $\mathcal{S}\subseteq R\text{-Simp}\cap \mathbb{T}_\sigma$ such that  $\xi_{\sigma\text{-nat}}(M)=\xi_{\sigma\text{-nat}}(\mathcal{S}).$ Now, since $L\in \xi_{\sigma\text{-nat}}(M),$ one has that there exists  $S\in\mathcal{S}$ such that $S$ embeds in $L.$ That is, $R$ is a left  $\sigma-$semiartinian ring.
\end{proof}

The following theorem generalizes \cite[Proposition 2.12]{on boolean}.

\begin{theorem}\label{sigmanat contenida en sigmaTORS}
Let  $\sigma\in R$-pr be exact and  costable. The following statements are equivalent::
\begin{enumerate}
    \item Each $S\in R\text{-Simp}\cap \mathbb{T}_\sigma$ is paraprojective.
    \item $\signat\subseteq \sigTors$.
    \item $\signat\subseteq \mathcal{L}_{/_\sigma}.$
\end{enumerate}
\end{theorem}
\begin{proof}
$(1)\Rightarrow(2)$ Let  $\mathcal{C}\in R\text{-}(\sigma\text{-nat})$. Since $\sigma\text{-}\left(R\text{-nat}\right)=\mathcal{L}_{\{\leq_\sigma,\oplus,\sigma(E( )),ext\}}$ and  $\sigma\text{-}\left(R\text{-TORS}\right)$\\ $=\mathcal{L}_{\{/_\sigma,\oplus, ext\}}$ it only remains to show that $\mathcal{C}\in\mathcal{L}_{/_\sigma}.$ To do so, let us consider $M\in\mathcal{C}$ and  $N\leq M$ with  $\sigma(N)\neq 0.$ Let $S\in\text{Simp}\cap\mathbb{T}_\sigma$ be a  submodule of $\sigma(N),$ which exists since $R$ is a left $\sigma$-semiartinian ring. If in the following  diagram $P$ is the pull-back of $M\twoheadrightarrow \sigma(N)$ and $S\rightarrowtail \sigma(N),$
\[\begin{tikzcd}
P \arrow[d, tail, dashed] \arrow[rr, two heads, dashed] & & S \arrow[d,tail] \\
M\arrow[r, two heads]& N \arrow[r, two heads] &\sigma(N).
\end{tikzcd}\]
Then $S$ is a quotient of $P,$ but by hypothesis $S$ is paraprojective, so $S$ is a simple  module  of  $\sigma$-torsion, which embeds in $P.$ Then $S$ is a   $\sigma$-torsion  simple module which embeds in  $M.$ The above shows that any $\sigma$-torsion simple submodule of     $\sigma(N)$ belongs to $\mathcal{C}$. As well, by  Proposition \ref{en sigmasemiar nat esta determinada por simples}, $\sigma(N)\in\mathcal{C}.$ Therefore $\mathcal{C}$ is  a $\sigma$-hereditary class.

$(2)\Rightarrow(3)$ It is clear.

$(3)\Rightarrow(1)$ Let $M$ be an $R$-module and $S\in\text{Simp}\cap\mathbb{T}_\sigma$ be a quotient of $M.$ By  hypothesis, one has that  $S\in\xi_{\sigma\text{-nat}}(M).$ It follows from Remark \ref{descripcion de la sigmanatural generada} that $S$ embeds in $M.$ That is, $S$ is paraprojective.
\end{proof}

\begin{corollary}\label{sigHH ent sigmanat contenida en sigmaTORS}
Let $\sigma\in R$-pr be exact and costable. If $R$ satisfies  condition $\sighh,$ then \[\signat\subseteq \sigTors.\]
\end{corollary}
\begin{proof}
It follows from statement \ref{SHH impli parapro} of Remark \ref{consecuen de SHH}, and from Theorem \ref{sigmanat contenida en sigmaTORS}.
\end{proof}

\begin{definition}
Let $\sigma\in R$-pr and let $_R M$ be an $R$-module.  We will say that a submodule $N$ of $M$ is $\sigma$-essential in $M$ if $N\cap L\neq 0$ for any  $L$ of $M$ with $\sigma(L)\neq 0$
\end{definition}

\begin{remark}\label{Zoc_sigma en semiart}
Let $\sigma\in R$-pr and let  $_R M$ be an  $R$-module. Let us denote $Soc_\sigma (M)$ the sum of the simple  $\sigma$-torsion submodules of \(M\).  If \(R\)  is a left \(\sigma\)-semiartinian ring then $Soc_\sigma (M)$ is essential in $\sigma (M),$ thus it is $\sigma$-essential in $_R M.$ Thus we have for a left $\sigma$-semiartinian ring that  $$\sigma(M)\neq 0 \Leftrightarrow Soc_\sigma (M)\neq0.$$
\end{remark}
\begin{theorem}
Let $\sigma\in R$-pr be exact and costable. If $R$ satisfies  condition $\sighh,$ then $\sigtors\subseteq \signat.$
\end{theorem}
\begin{proof} Let $ \mathcal{C}\in \sigtors.$   To see that  $\mathcal{C}\in\signat $ recall that $\signat =\mathcal{L}_{\leq_\sigma,\oplus,\sigma E}$  and that $\sigtors = \mathcal{L}_{\leq_\sigma,\oplus, /_\sigma , ext}$. Thus, it suffices to  show that if $M\in \mathcal{C}$, then $\sigma E(M)\in\mathcal{C}.$ 
Let us first note that for any modulo $_RN$ one has that $$Soc_\sigma (N)=Soc_\sigma (E(N)),$$ since $N$ is an essential submodule of $E(N) $.
On the other hand, as $\sigma$ is an exact and stable preradical, and   in particular it is  left exact  and stable,  we have that $$\sigma(E(M))=\sigma(E(M)).$$ From Lemma 1.42 and the equations one has that
\begin{align*}
    M\in\mathcal{C} & \Longleftrightarrow Soc_{\sigma}(M)\in\mathcal{C}\\
    & \Longleftrightarrow Soc_{\sigma}(E\left(M\right))\in\mathcal{C}\\
    & \Longleftrightarrow Soc_{\sigma}(\sigma\left(E\left(M\right)\right))\in C\\
    & \Longleftrightarrow\sigma\left(E\left(M\right)\right)\in\mathcal{C}.
\end{align*}
Thus, if $M\in\mathcal{C},$ then $E(M)\in\mathcal{C}.$
\end{proof}

\begin{definition} Let $\sigma\in R-pr$ exact an stable and assume that $R$ has $\sighh.$  We are going to define for each ordinal $\alpha$, a preradical $Soc_{\sigma}^{\alpha}$. Let us define $Soc_{\sigma}^{0}$ as the zero preradical $\text{\ensuremath{\underbar{0}}}.$  Now, let us define $Soc_{\sigma}^{\alpha+1}\left(M\right)/Soc_{\sigma}^{\alpha}\left(M\right)=Soc_{\sigma}\left(M/Soc_{\sigma}^{\alpha}\left(M\right)\right),$ as in the diagram
\begin{equation}\label{SLoewy}
\begin{tikzcd}
0 \arrow[r, ] & Soc_{\sigma}^{\alpha}\left(M\right))\arrow[r,  ] &  
M\arrow[r,  ] & M/Soc_{\sigma}^{\alpha}\left(M\right)\arrow[r, ]&0\\0\arrow[r, ]& Soc_{\sigma}^{\alpha}\left(M\right))\arrow[r,  ]\arrow[u, equal]&Soc_{\sigma}^{\alpha+1}\left(M\right) \arrow[r, ]\arrow[u, hook]&Soc_{\sigma}\left(M/Soc_{\sigma}^{\alpha}\left(M\right)\right)\arrow[r, ]\arrow[u, hook]&0.
\end{tikzcd}
\end{equation} Let us define $Soc_{\sigma}^{\alpha}\left(M\right)=\sum_{\beta<\alpha}Soc_{\sigma}^{\beta}\left(M\right)$ if   $\alpha$ is a limit ordinal.

\end{definition}   
\begin{lemma}  Let $\sigma\text{\ensuremath\in R-pr}$ be exact and stable and let  $\mathcal{C}$ be a $\sigma$-hereditary torsion class. If $R$ has the $\sighh$ condition, then $M\in\mathcal{C}$ if and only if each of its simple $\sigma$-torsion subquotients belong to $\mathcal{C}.$

\end{lemma}
\begin{proof}
Let us first assume that $M\in\mathcal{C}.$ If $\sigma(M)\neq0$ and $S\in R-simp\cap\mathbb{T_{\sigma}}$ is a subquotient of $\sigma\left(M\right)$, then by Theorem 1,10, $S$ embeds in $\sigma\left(M\right).$ As $\mathcal{C}$ is a $\sigma$-hereditary class, one has that $S\in\mathcal{C\mathrm{.}}$ 

On the other hand, if $\sigma(M)=0,$ then $M\in\mathcal{C},$ because $\mathbb{F_{\sigma}}\subseteq\mathcal{C},$ for any $\sigma$-hereditary torsion class $\mathcal{C}.$

Now suppose that any simple $\sigma$-torsion subquotient of $\sigma(M)$ belongs to $\mathcal{C}$ and let us show that $M$ in $\mathcal{C}.$ 

 We can assume that $\sigma(M)\neq0.$ Since $\sigma$ is exact and stable, then it centrally splits, so that $M=\sigma(M)\oplus M',$ with $M'\in\mathbb{F}_{\sigma}.$ Thus, it suffices to show that $\sigma(M)\in\mathcal{C}.$ As  $R$ has $ \sighh,$ then by Corollary \ref{SHH Smax y Ssemi}, we have that $R$ is left $\sigma$-semiartinian. Then by Remark \ref{Zoc_sigma en semiart}, $0\neq Soc_{\sigma}\leqslant_{ess}\sigma(M).$ Note also that $Soc_{\sigma}(M)\in\mathcal{C},$ since $\mathcal{C}$ is a $\sigma$-hereditary torsion class. Therefore, if it happens that $\sigma(M)=Soc_{\sigma}(M),$ then $\sigma(M)$ belongs to $\mathcal{C}.$ In this case we are done. Suppose then that $\sigma(M)/Soc_{\sigma}(M)\neq0.$ Note that we can consider two cases: $\sigma(M)/Soc_{\sigma}(M)\text{\ensuremath{\in\mathbb{F_{\sigma}}}}$ and $\sigma(M)/Soc_{\sigma}(M)\text{\ensuremath{\notin\mathbb{F_{\sigma}}}}.$ 

In the first case, $\sigma(M)/Soc_{\sigma}(M)\in\mathcal{C}$ because $\mathbb{F}_{\sigma}\subseteq\mathcal{C}.$ Since $\mathcal{C}$ is closed under extensions, the exact sequence $0\rightarrow Soc_{\sigma}(M)\rightarrow\sigma(M)\rightarrow\sigma(M)/Soc_{\sigma}(M)\rightarrow0$ implies that $\sigma(M)\in\mathcal{C}.$

In the second case, we are going to show by induction that for any ordinal $\alpha,$ the module $Soc_{\sigma}^{\alpha}\left(\sigma\left(M\right)\right)$ belongs to $\mathcal{C}.$

It is clear that $$\text{\ensuremath{\underbar{0}}}\left(\sigma\left(M\right)\right)=Soc_{\sigma}^{0}\left(\sigma\left(M\right)\right)=_{R}0$$ belongs to $\mathcal{C}.$ Now suppose that $Soc_{\sigma}^{\alpha}\left(\sigma\left(M\right)\right)$ belongs to $\mathcal{C}$,  then the exact sequence $$0\rightarrow Soc_{\sigma}^{\alpha}\left(\sigma\left(M\right)\right)\rightarrow Soc_{\sigma}^{\alpha+1}\left(\sigma\left(M\right)\right)\rightarrow Soc_{\sigma}\left(\sigma\left(M\right)/Soc_{\sigma}^{\alpha}\left(\sigma\left(M\right)\right)\right)\rightarrow0$$ has the ends in $C$, because any simple submodule of $Soc_{\sigma}\left(\sigma\left(M\right)/Soc_{\sigma}^{\alpha}\left(\sigma\left(M\right)\right)\right)$ is a $\sigma$-torsion simple subquotient of $\sigma\left(M\right)$, which is in $\mathcal{C},$ by hypothesis.

So, $Soc_{\sigma}\left(\sigma\left(M\right)/Soc_{\sigma}^{\alpha}\left(\sigma\left(M\right)\right)\right)$ belongs to $\mathcal{C},$ because $\mathcal{C}$ is closed under direct sums and quotients of $\sigma$-torsion modules.
 If $\alpha$ is a limit ordinal, suppose that $Soc_{\sigma}^{\beta}\left(\sigma\left(M\right)\right)\in\mathcal{C}$ for each $\beta<\alpha,$ then $Soc_{\sigma}^{\alpha}\left(\sigma\left(M\right)\right)=\sum_{\beta<\alpha}Soc_{\sigma}^{\beta}\left(\sigma\left(M\right)\right)$ is also in $\mathcal{C},$because this sum is a quotient of $\oplus_{\beta<\alpha}Soc_{\sigma}^{\beta}\left(\sigma\left(M\right)\right).$  As $\left\{ Soc_{\sigma}^{\alpha}\left(\sigma\left(M\right)\right)\right\}$  is an ascending chain of submodules of $\sigma\left(M\right)$ which can not be all distinct, thus there exists  a first ordinal $\gamma$ such that $ Soc_{\sigma}^{\gamma}\left(\sigma\left(M\right)\right)= Soc_{\sigma}^{\gamma+1}\left(\sigma\left(M\right)\right)\} \} .$ This means that $0=Soc_{\sigma}^{\gamma+1}\left(\sigma\left(M\right)\right)/Soc_{\sigma}^{\gamma}\left(\sigma\left(M\right)\right)=Soc_{\sigma}\mathcal{\left(\mathit{M/Soc_{\sigma}^{\gamma}\left(\sigma\left(M\right)\right)}\right)}.$ This can only occur if $M/Soc_{\sigma}^{\gamma}\left(\sigma\left(M\right)\right)=0.$ Therefore $M=Soc_{\sigma}^{\gamma}\left(\sigma\left(M\right)\right),$ which belongs to $\mathcal{C\mathit{.}}$
It is clear that $\text{\ensuremath{\underbar{0}}}\left(\sigma\left(M\right)\right)=Soc_{\sigma}^{0}\left(\sigma\left(M\right)\right)=_{R}0$ belongs to $\mathcal{C}.$ Now suppose that $Soc_{\sigma}^{\alpha}\left(\sigma\left(M\right)\right),$ belongs to $\mathcal{C}$,  then the exact sequence $ 0\rightarrow Soc_{\sigma}^{\alpha}\left(\sigma\left(M\right)\right)\rightarrow Soc_{\sigma}^{\alpha+1}\left(\sigma\left(M\right)\right)\rightarrow Soc_{\sigma}\left(\sigma\left(M\right)/Soc_{\sigma}^{\alpha}\left(\sigma\left(M\right)\right)\right)\rightarrow0$ has the ends in $\mathcal{C}$, because any simple submodule of $Soc_{\sigma}\left(\sigma\left(M\right)/Soc_{\sigma}^{\alpha}\left(\sigma\left(M\right)\right)\right)$ is a $\sigma$-torsion simple subquotient of $\left(M\right),$ which is in $\mathcal{C},$ by hypothesis. So $Soc_{\sigma}\left(\sigma\left(M\right)/Soc_{\sigma}^{\alpha}\left(\sigma\left(M\right)\right)\right)$ belongs to $\mathcal{C}$, because $\mathcal{C}$ is closed under direct sums and quotients of $\sigma$-torsion modules.
If $\alpha$ is a limit ordinal, suppose that $Soc_{\sigma}^{\beta}\left(\sigma\left(M\right)\right)\in\mathcal{C}$ for each $\beta<\alpha,$ then $Soc_{\sigma}^{\alpha}\left(\sigma\left(M\right)\right)=\sum_{\beta<\alpha}Soc_{\sigma}^{\beta}\left(\sigma\left(M\right)\right)$ is also in $\mathcal{C},$ because the sum is a quotient of $\oplus_{\beta<\alpha}Soc_{\sigma}^{\beta}\left(\sigma\left(M\right)\right).$

As $\left\{ Soc_{\sigma}^{\alpha}\left(\sigma\left(M\right)\right)\right\}$  is an ascending chain of submodules of $\sigma\left(M\right)$ which can not all be  distinct  there exists a first ordinal $\gamma$  such that $ Soc_{\sigma}^{\gamma}\left(\sigma\left(M\right)\right)= Soc_{\sigma}^{\gamma+1}\left(\sigma\left(M\right)\right)  .$ This means that \[0=Soc_{\sigma}^{\gamma+1}\left(\sigma\left(M\right)\right)/Soc_{\sigma}^{\gamma}\left(\sigma\left(M\right)\right)=Soc_{\sigma}\mathcal{\left(\mathit{M/Soc_{\sigma}^{\gamma}\left(\sigma\left(M\right)\right)}\right)}.\] This can only occur if $M/Soc_{\sigma}^{\gamma}\left(\sigma\left(M\right)\right)=0.$ Therefore $M=Soc_{\sigma}^{\gamma}\left(\sigma\left(M\right)\right),$ which belongs to $\mathcal{C\mathit{.}}$
\end{proof}
\begin{theorem}\label{sigHH ent signat=sigconat}
Let $\sigma\in R$-pr be an exact costable preradical. If  $R$ satisfies condition $\sighh,$ then
\[\sigma\text{-}\left(R\text{-nat}\right)=\sigma\text{-}\left(R\text{-conat}\right).\]
\end{theorem}
\begin{proof}
$\subseteq]$ Let us take a   $\sigma$-natural class  $\mathcal{C}$,  i.e.  a module class satisfying  $\left(\sigma\text{-}N\right).$ We will show that   $\mathcal{C}$ satisfies $\left(\sigma\text{-}CN\right),$  that is, we are going to show that if  for an $R$-module $_{R}M$ and for   each non zero  epimorphism   $M\twoheadrightarrow L,$ with  $\sigma(L)\neq 0,$  there exist $C\in\mathcal{C}$ and  $N\in R$-Mod, with $\sigma(N)\neq 0,$ and two epimorphisms  $L\twoheadrightarrow N\twoheadleftarrow C,$ then  $M\in\mathcal{C}.$

Let  $f:K\rightarrowtail M$ be  a monomorphism with $\sigma(K)\neq 0.$ Since $R$ satisfies  condition $\sighh,$ then by  Corollary \ref{SHH Smax y Ssemi}, one has that $R$ is a left  $\sigma$-semiartinian ring, so that there exists $S\in R\text{-}Simp\cap\mathbb{T}_\sigma$ such that  $S$ embeds in $K$ and hence embeds in $M.$ By Remark \ref{consecuen de SHH}, one has that  $S$ is parainjective, so $S$ is a quotient of $M.$ Now, by hypothesis, we have that there exists $C\in\mathcal{C}$ such that $S$ is a quotient of $C.$ 
By Remark \ref{consecuen de SHH}, we have that  $S$ is paraprojective, so $S$ embeds in $C.$ Note that we find ourselves in the following situation
\begin{equation}\label{Snat}
\begin{tikzcd}
K \arrow[r, tail] & M &  \\
& S\arrow[u, tail] \arrow[r, tail]& C.
\end{tikzcd}
\end{equation}
Since $\mathcal{C}$ satisfies the condition $\left(\sigma\text{-}N\right),$ from \eqref{Snat}, it follows that $M\in\mathcal{C}.$
That is to say, $\mathcal{C}\in\sigma\text{-}\left(R\text{-conat}\right).$
\bigskip

\noindent $\supseteq]$ Let $\mathcal{C}$ be a $\sigma$-conatural class. We will show that $\mathcal{C}$ satisfies the $\left(\sigma\text{-}N\right)$ condition, i.e, we will show that if $M$ is a $R$-module such that for any monomorphism  $L\rightarrowtail M,$ with $\sigma(L)\neq 0$ there exist $C\in\mathcal{C}$ and
 $N\in R$-Mod, with $\sigma(N)\neq 0,$ which satisfy that $L\leftarrowtail N\rightarrowtail C,$ then $M\in\mathcal{C}.$

Consider then an   epimorphism $f:M\twoheadrightarrow K$, with $\sigma(K)\neq 0.$ Since $R$ satisfies  condition $\sighh,$ by  Corollary \ref{SHH Smax y Ssemi}, $R$ is a left $\sigma$-max ring, so there exists $S\in R\text{-}Simp\cap\mathbb{T}_\sigma$ such that $S$ is a quotient of $K$ and hence a quotient of $M.$ Now, by Remark \ref{consecuen de SHH}, one has that $S$ is paraprojective, so  $S$ embeds in $M.$ Then, by hypothesis, we have that there exists $C\in\mathcal{C}$ such that  $S$ embeds in $C.$ 
Note that, from Remark \ref{consecuen de SHH}, one has that  $S$ is parainjective, so $S$ is a quotient of $C.$ Thus, we have the following situation
\begin{equation}\label{Sconat}
\begin{tikzcd}
M \arrow[r, two heads] & K\arrow[d, two heads] &  \\
& S & C.\arrow[l, two heads]
\end{tikzcd}
\end{equation}
Now, since $\mathcal{C}\in\sigma\text{-}\left(R\text{-conat}\right),$ $\mathcal{C}$ satisfies condition $\left(\sigma\text{-}CN\right),$ so, from \eqref{Sconat}, one has that $M\in\mathcal{C}.$
Thus, $\mathcal{C}\in\sigma\text{-}\left(R\text{-nat}\right).$
\end{proof}

\begin{theorem}
Let $\sigma\in R$-pr be exact and costable. If $\sigma\text{-}\left(R\text{-nat}\right)=\sigma\text{-}\left(R\text{-conat}\right),$ then any $S\in R\text{-}Simp \cap\mathbb{T}_\sigma$ is parainjective and paraprojective.
\end{theorem}
\begin{proof}
Let $M$ be a $R$-module and $S\in R\text{-}Simp \cap\mathbb{T}_\sigma$. If $S$ embeds in $M,$ since $S=\sigma(S),$ then $S\in\xi_{\sigma\text{-}nat}(M).$ Now, since  $\sigma\text{-}\left(R\text{-nat}\right)=\sigma\text{-}\left(R\text{-conat}\right),$ it follows that $\xi_{\sigma\text{-}nat}(M)=\xi_{\sigma\text{-}conat}(M).$ Thus, $S\in\xi_{\sigma\text{-}conat}(M)$, then, by \ref{descripcion de la sigmaconatural generada}, $S$ is a quotient of $M,$ i.e., $S$ is parainjective.

On the other hand, if  $S$ is a quotient of $M,$ then $S\in\xi_{\sigma\text{-}conat}(M),$ since $S=\sigma(S).$ It follows from the fact that  $\sigma\text{-}\left(R\text{-nat}\right)=\sigma\text{-}\left(R\text{-conat}\right),$ that $S\in\xi_{\sigma\text{-}nat}(M).$ By \ref{descripcion de la sigmanatural generada}, one has that $S$ embeds in $M,$ i.e., $S$ is paraprojective. If
\end{proof}

\begin{remark}\label{sigmaR-módulos}
Sea $\sigma\in R$-pr. If $\sigma$ centrally splits , then there  exists a central idempotent $e\in R,$ such that  $\sigma(\_)=e\cdot\_.$ Now, as  $R=Re\oplus R(1-e),$ it follows that $\sigma(R)\cong R/(1-e)R,$ so a $\sigma(R)$-module is an $R$-module annulled by $1-e$. As well, if  $M$ i s a $\sigma(R)$-module, then $M=\sigma(M).$ 
\end{remark}

\begin{lemma}\label{HH sii sigHH}
Let $\sigma\in R$-pr. If $\sigma$ centrally splits, then the following statements are equivalent:
\begin{enumerate}

\item $R$ satisfies the  condition $\sighh$.
\item      $\sigma(R)$ satisfies the condition $(HH).$ 
\end{enumerate}
\end{lemma}
\begin{proof}
$(1)\Rightarrow(2)$ Let $M,N$ be two $\sigma(R)$-modules and let  $0\neq f\in Hom_{\sigma(R)}(M,N).$ Since  $M$ and $N$ are also $R$-modules, it follows that  $0\neq f\in Hom_R(M,N)$. Now, by Remark \ref{sigmaR-módulos}, it follows that $M=\sigma(M)$, so that   $0\neq f\in Hom_R(\sigma(M),N).$ Since  $R$ satisfies the condition $\sighh$,  it follows that there exists $0\neq g\in Hom_R(N,\sigma(M))$, whence  $0\neq g\in Hom_{\sigma(R)}(N,M).$That is , $\sigma(R)$ satisfies the condition $\sighh$.

$(2)\Rightarrow(1)$ Let $M,N$ be two  $R$-modules such that  $\sigma(M),\sigma(N)\neq 0$ and  $0\neq f\in Hom_R(\sigma(M),N).$ Note that $\sigma(M)$ and  $\sigma(N)$ are two $\sigma(R)$-modules. Now, since $\sigma$ is a preradical, then $f$ can be corestricted to $\sigma(N),$ i.e., $f\lceil:\sigma(M)\to\sigma(N)$ is a  $\sigma(R)$ nonzero morphism. It follows that there exists $0\neq g\in Hom_{\sigma(R)}(\sigma(N),\sigma(M))$, since $\sigma(R)$ satisfies the  condition $(HH)$. Thus, we find ourselves in the following situation. 
\[\begin{tikzcd}
\sigma(M) \ar[dd, equal] \ar[rr, "f\lceil"] & & \sigma(N) \ar[lldd, "h\neq 0"] \ar[dd, hook] \\
& \\
\sigma(M) & & N=\sigma(N)\oplus(1-e)N. \ar[ll, dashed, "h\oplus\overline{0}"]
\end{tikzcd}\]
Then, if $\overline{0}$ denotes the zero morphism, then we have that  $h\oplus\overline{0}$ is a nonzero morphism from $N$ to $\sigma(M)$. Therefore, $R$ satisfies the  condition $\sighh$.
\end{proof}

The next theorem  generalizes \cite[Theorem 2.18]{On retractability and its}.

\begin{theorem}
Let $R$ be a ring and $\sigma\in R\text{-}pr$ be exact and costable. The following statements are equivalent:
\begin{enumerate}
    \item $R$ satisfies condition $\sighh.$
    \item $R$ is a  $\sigmod$-retractable ring and  each  $S\in  R\text{-}Simp\cap\mathbb{T}_\sigma$ is paraprojective.
    \item Each $S\in  R\text{-}Simp\cap\mathbb{T}_\sigma$ is paraprojective and parainjective.
    \item $ \sigTors= \sigconat=\sigtors=\signat.$
    \item  $\signat=\sigconat.$
    \item $R=\sigma(R)\times R_2,$ where $\sigma(R)$ is a ring that satisfies the $(HH)$ condition.
    \item $R=\sigma(R)\times R_2,$ where $\sigma(R)$ is isomorphic to a full matrix ring over a local left and right perfect ring.
\end{enumerate}
\end{theorem}
\begin{proof}
$(1)\Leftrightarrow (2)\Leftrightarrow (3)$ is Theorem \ref{sigHH sii pariny y parapro}.
$(1)\Rightarrow (5)$ is Theorem \ref{sigHH ent signat=sigconat}.
$(5)\Rightarrow (1)$ It follows from Theorem \ref{sigHH sii pariny y parapro}.
$(1)\Rightarrow (4)$ is  Theorem \ref{sigHH ent sigTors=sigtors=sigconat}. $(4) \Rightarrow (5)$ This is clear. $(1)\Leftrightarrow(6)$ This follows from Lemma \ref{HH sii sigHH} and because $\sigma$ centrally splits. $(6)\Leftrightarrow (7)$ It follows from \cite[Theorem 2.18]{On retractability and its}.
\end{proof}


\providecommand{\bysame}{\leavevmode\hbox
to3em{\hrulefill}\thinspace}


\end{document}